\documentclass[11pt, reqno]{amsart}

\usepackage[OT2, T1]{fontenc}
\usepackage{url}
\usepackage{amsmath}
\usepackage{array}
\usepackage{graphicx}
\usepackage{amsfonts}
\usepackage{amssymb}
\usepackage{amsthm}
\usepackage{hyperref}
\usepackage{colonequals}	
\usepackage{color}
\usepackage{mathabx}		
\usepackage{mathrsfs}
\usepackage{amscd}
\usepackage[all,cmtip]{xy}	
\usepackage{sseq}			
\usepackage{verbatim}
\usepackage{parskip}
\usepackage[in]{fullpage}
\usepackage{mathrsfs} 			
\usepackage{bm} 				
\usepackage{cleveref}
\usepackage{enumitem}
\usepackage[alphabetic,lite]{amsrefs} 	

\DeclareSymbolFont{cyrletters}{OT2}{wncyr}{m}{n}
\DeclareMathSymbol{\Sha}{\mathalpha}{cyrletters}{"58}

\usepackage{color}

\newcommand{\gA}{\alpha}

\newcommand{\gL}{\lambda}

\newcommand{\bA}{\mathbb{A}}
\newcommand{\bC}{\mathbb{C}}

\newcommand{\bF}{\mathbb{F}}
\newcommand{\bG}{\mathbb{G}}

\newcommand{\bQ}{\mathbb{Q}}
\newcommand{\bR}{\mathbb{R}}

\newcommand{\bZ}{\mathbb{Z}}

\newcommand{\cE}{\mathcal{E}}
\newcommand{\cF}{\mathcal{F}}

\newcommand{\cO}{\mathcal{O}}

\newcommand{\cW}{\mathcal{W}}
\newcommand{\cX}{\mathcal{X}}

\newcommand{\fm}{\mathfrak{m}}

\newcommand{\fp}{\mathfrak{p}}

\newcommand{\ra}{\rightarrow}

\newcommand{\xra}{\xrightarrow}

\newcommand{\hra}{\hookrightarrow}

\newcommand{\I}{^{\infty}}
\newcommand{\wt}{\widetilde}
\newcommand{\wh}{\widehat}
\newcommand{\eps}{\epsilon}

\newcommand{\pr}{^{\prime}}
\newcommand{\ce}{\colonequals}
\newcommand{\ov}{\overline}

\renewcommand{\b}{\textbf}

\newcommand{\tensor}{\otimes} 		
\newcommand{\isomto}{\overset{\sim}{\longrightarrow}}

\renewcommand{\i}{^{-1}}
\renewcommand{\div}{_{\text{div}}}		
\newcommand{\leftexp}[2]{{\vphantom{#2}}^{#1}{#2}}

\providecommand{\abs}[1]{\left\lvert#1\right\rvert}

\providecommand{\p}[1]{\left(#1\right)}

\providecommand{\f}[2]{\frac{#1}{#2}}
\DeclareMathOperator{\Ker}{Ker}			
\DeclareMathOperator{\coker}{coker}		
\DeclareMathOperator{\Spec}{Spec}		
\DeclareMathOperator{\Hom}{Hom}			
\DeclareMathOperator{\Char}{char}		
\DeclareMathOperator{\Frob}{Frob}		
\DeclareMathOperator{\Gal}{Gal}	
\DeclareMathOperator{\tr}{tr}		
\DeclareMathOperator{\ord}{ord}	
\DeclareMathOperator{\ab}{ab}		
\DeclareMathOperator{\Ind}{Ind}		
\DeclareMathOperator{\End}{End}		
\DeclareMathOperator{\Aut}{Aut}		
\DeclareMathOperator{\rk}{rk}		
\DeclareMathOperator{\Sel}{Sel}		
\DeclareMathOperator{\Art}{Art}		
\DeclareMathOperator{\Lie}{Lie}		
\DeclareMathOperator{\length}{length}		
\newcommand{\ba}{\begin{aligned}}
\newcommand{\ea}{\end{aligned}}
\newcommand{\be}{\begin{equation}}
\newcommand{\ee}{\end{equation}}
\newcommand{\pf}{\begin{proof}}
\newcommand{\bpf}{\begin{proof}}
\newcommand{\epf}{\end{proof}}
\newcommand{\bthm}{\begin{thm}}
\newcommand{\ethm}{\end{thm}}
\newcommand{\bprop}{\begin{prop}}
\newcommand{\eprop}{\end{prop}}
\newcommand{\bcor}{\begin{cor}}
\newcommand{\ecor}{\end{cor}}
\newcommand{\brem}{\begin{rem}}
\newcommand{\erem}{\end{rem}}
\newcommand{\blem}{\begin{lemma}}
\newcommand{\elem}{\end{lemma}}
\newcommand{\bconj}{\begin{conj}}
\newcommand{\econj}{\end{conj}}
\newcommand{\benum}{\begin{enumerate}}
\newcommand{\eenum}{\end{enumerate}}
\newcommand{\bc}{}
\newcommand{\lab}{\label}

\theoremstyle{plain}
\newtheorem{thm}[subsection]{Theorem}
\Crefname{thm}{Theorem}{Theorems}
\newtheorem{rethm}{Theorem}
\Crefname{rethm}{Theorem}{Theorem}
\newtheorem{prop}[subsection]{Proposition}
\Crefname{prop}{Proposition}{Propositions}

\Crefname{Problem}{Problem}{Problems}
\newtheorem{conj}[subsection]{Conjecture}
\Crefname{conj}{Conjecture}{Conjectures}
\newtheorem{cor}[subsection]{Corollary}
\Crefname{cor}{Corollary}{Corollaries}

\newtheorem{lemma}[subsection]{Lemma}

\theoremstyle{remark}
\newtheorem{claim}[subsubsection]{Claim}
\Crefname{claim}{Claim}{Claims}

\theoremstyle{definition}

\theoremstyle{definition}

\newtheorem{rem}[subsection]{Remark}
\Crefname{rem}{Remark}{Remarks}

\newtheoremstyle{subsection-tweak}
   {11pt}
   {3pt}%
   {}
   {}%
   {\bfseries}
   {}%
   {.5em}
   {\thmnumber{\@{#1}{}\@{#2}.}%
    \thmnote{~{\bfseries#3.}}}    
    
\newtheorem{innercustomconj}{Conjecture}
\newenvironment{customconj}[1]
  {\renewcommand\theinnercustomconj{#1}\innercustomconj}
  {\endinnercustomconj}

\Crefname{innercustomconj}{Conjecture}{Conjecture}

\theoremstyle{subsection-tweak}
\newtheorem{pp}[subsection]{}
\newcommand{\bpp}{\begin{pp}}
\newcommand{\epp}{\end{pp}}

\begin{document}
\author{K\k{e}stutis \v{C}esnavi\v{c}ius}
\title{The $p$-parity conjecture for elliptic curves with a $p$-isogeny}
\date{\today}
\address{Department of Mathematics, Massachusetts Institute of Technology, Cambridge, MA 02139, USA}
\email{kestutis@math.mit.edu}
\urladdr{http://math.mit.edu/~kestutis/}

\begin{abstract} For an elliptic curve $E$ over a number field $K$, one consequence of the Birch and Swinnerton-Dyer conjecture is the parity conjecture: the global root number matches the parity of the Mordell--Weil rank. Assuming finiteness of $\Sha(E/K)[p^\infty]$ for a prime $p$ this is equivalent to the $p$-parity conjecture: the global root number matches the parity of the $\bZ_p$-corank of the $p\I$-Selmer group. We complete the proof of the $p$-parity conjecture for elliptic curves that have a $p$-isogeny for $p > 3$ (the cases $p \le 3$ were known). T. and V. Dokchitser have showed this in the case when $E$ has semistable reduction at all places above $p$ by establishing respective cases of a conjectural formula for the local root number. We remove the restrictions on reduction types by proving their formula in the remaining cases. We apply our result to show that the $p$-parity conjecture holds for every $E$ with complex multiplication defined over $K$. Consequently, if for such an elliptic curve $\Sha(E/K)[p^\infty]$ is infinite, it must contain $(\bQ_p/\bZ_p)^2$. \end{abstract}

\maketitle

\section{Introduction}\lab{S0}

If $E$ is an elliptic curve defined over a number field $K$, its completed $L$-series is conjectured to have a holomorphic continuation $\Lambda(E/K, s)$ to the whole complex plane, and to satisfy a functional equation
\be\lab{A} \Lambda(E/K, 2 - s) = w(E/K)\Lambda(E/K, s).\ee
Here $w(E/K) \in \{ \pm 1\}$ is the \emph{global root number} of $E/K$. It can be given a definition independent of \eqref{A} as the product 
\be\lab{B} w(E/K) = \prod_{v\text{ place of } K} w(E/K_v)\ee
of \emph{local root numbers} $w(E/K_v) \in \{ \pm 1\}$ (here $K_v$ is the completion of $K$ at $v$), with $w(E/K_v)$ defined as the root number of the Weil--Deligne representation associated to $E/K_v$ if $v\nmid \infty$, and $w(E/K_v) = -1$ if $v\mid \infty$ (cf., for instance, \cite{Roh94}).

Granting holomorphic continuation of $\Lambda(E/K, s)$, the Birch and Swinnerton-Dyer conjecture (BSD) predicts that
\be\lab{C} \ord_{s = 1} \Lambda(E/K, s) = \rk E(K),\ee
where $\rk E(K) \ce \dim_{\bQ} E(K) \tensor \bQ$ is the \emph{Mordell--Weil rank} of $E/K$. Combining \eqref{A} and \eqref{C} one gets

\bconj[Parity conjecture]\lab{D} $(-1)^{\rk E(K)} = w(E/K)$.\econj

The parity conjecture is more approachable than BSD, and Tim and Vladimir Dokchitser have showed \cite{DD11}*{Theorem 1.2} that it holds if one assumes that the $2$- and $3$-primary parts of the Shafarevich--Tate group $\Sha(E/K(E[2]))$ are finite ($K(E[2])$ is the smallest extension of $K$ over which the $2$-torsion of $E$ is rational). 

If one hopes for unconditional results, one is led to consider the $p\I$-Selmer rank $\rk_p(E/K)$ instead of $\rk E(K)$ for each prime $p$. To define it one takes the exact sequence
\be\lab{E} 
0 \ra E(K) \tensor \bQ_p/\bZ_p \ra \varinjlim \Sel_{p^n}(E/K) \ra \Sha(E/K)[p\I] \ra 0  
\ee
and lets $\rk_p(E/K) \ce \dim_{\bQ_p} \Hom_{\bZ_p} (\varinjlim \Sel_{p^n}(E/K), \bQ_p/\bZ_p) \tensor_{\bZ_p} \bQ_p$ be the $\bZ_p$-corank of the $p$-primary torsion abelian group $\varinjlim \Sel_{p^n}(E/K)$. From \eqref{E} one gets 
\be\lab{rk-sum}
 \rk_p(E/K) = \rk E(K) + r_p(E/K),
 \ee
where $r_p(E/K)$ is the $\bZ_p$-corank of $\Sha(E/K)[p\I]$ (equivalently, $r_p(E/K) = $ the number of copies of $\bQ_p/\bZ_p$ in $\Sha(E/K)[p\I]$). Since $\Sha(E/K)$ is conjectured to be finite \cite{Tat74}*{Conjecture 1} (to the effect that $r_p = 0$), \Cref{D} leads to

\bconj[$p$-parity conjecture] $(-1)^{\rk_p(E/K)} = w(E/K).$ \econj

The $p$-parity conjecture is known if $K = \bQ$ thanks to the work of Nekov\'{a}\v{r} \cite{Nek06}*{\S0.17}, Kim \cite{Kim07}, and T. and V. Dokchitser \cite{DD10}*{Theorem 1.4}, and also if $K$ is totally real excluding some cases of potential complex multiplication \cite{Nek09}, \cite{Nek12}*{Theorem A}, \cite{Nek14}*{5.12}. Over arbitrary $K$ and for arbitrary $p$ the following theorem of T. and V. Dokchitser is the most general result currently known (one can weaken the assumptions on reduction types above $p$ somewhat, see \Cref{L}).

\bthm[\cite{DD08}*{Theorem 2}, \cite{DD11}*{Corollary 5.8}]\lab{F} The $p$-parity conjecture holds for $E/K$ provided that $E$ has a $p$-isogeny (defined over $K$) and either $p \le 3$ or $E$ has semistable reduction at all places of $K$ above $p$.\ethm

The main goal of this paper is to remove the semistable hypothesis in \Cref{F}. Namely, we complete the proof of the following result.

\bthm[\S\ref{S18}, \S\ref{p3}, \S\ref{p2}, and \Cref{DONE}]\lab{G} The $p$-parity conjecture holds for $E/K$ provided that $E$ has a $p$-isogeny. \ethm

In fact, both the global root number and the parity of the $p\I$-Selmer rank do not change in an odd degree Galois extension \cite{DD09}*{Proposition A.2 (3)}, so this gives a slightly stronger

\renewcommand{\therethm}{\arabic{section}.\arabic{subsection}$^\prime$}
\begin{rethm}\lab{H} The $p$-parity conjecture holds for $E/K$ provided that $E$ acquires a $p$-isogeny over an odd degree Galois extension of $K$. \end{rethm}

Since $\rk E(K) = \rk_p(E/K)$ is equivalent to finiteness of $\Sha(E/K)[p\I]$, from \Cref{H} we get

\bcor The parity conjecture holds for $E/K$ provided that $E$ acquires a $p$-isogeny over an odd degree Galois extension of $K$ and $\Sha(E/K)[p\I]$ is finite.\ecor

\bigskip
\large{{\it{Implications for elliptic curves with complex multiplication}}}

\normalsize

I thank Karl Rubin for pointing out to me that one corollary of \Cref{G} is

\bthm[\Cref{CM-A}] \lab{CM-Aa} Let $E$ be an elliptic curve defined over a number field $K$, and suppose that $\End_K E \neq \bZ$, i.e., that $E$ has complex multiplication defined over $K$. Then the $p$-parity conjecture holds for $E/K$ for every prime $p$. \ethm

\brem If $E$ is an elliptic curve over a totally real field $L$ with complex multiplication (necessarily defined over a non-trivial extension of $L$), the $p$-parity conjecture for $E/L$ has been established by Nekov{\'a}{\v{r}} \cite{Nek12}*{Theorem A} in the cases where $2 \nmid [L : \bQ]$ or $p$ splits in $(\End_{\ov{L}} E) \tensor \bQ$.  \erem

If $E/K$ is as in \Cref{CM-Aa}, then $F \ce (\End_K E) \tensor \bQ$ is an imaginary quadratic field and $E(K) \tensor \bQ$ is an $F$-vector space. Consequently, $\rk E(K) = \dim_\bQ E(K) \tensor \bQ$ is even. This is half of

\bthm[\Cref{CM-B}]\lab{CM-Ba} If $E$ has complex multiplication defined over $K$, then the parity conjecture holds for $E/K$. More precisely, $\rk E(K)$ is even and $w(E/K) = 1$. \ethm

Since $w(E/K) = 1$, \Cref{CM-Aa} tells us that $\rk_p(E/K)$ is even. Therefore, from \eqref{rk-sum} we obtain that $r_p(E/K)$ is even as well. Concerning the conjectural finiteness of $\Sha(E/K)$, this gives

\bthm If $E$ has complex multiplication defined over $K$ and $\Sha(E/K)[p\I]$ is infinite, then $\Sha(E/K)[p\I]$ contains $(\bQ_p/\bZ_p)^2$. \ethm

If $\Sha(E/K)_{\div}$ denotes the divisible part of the Shafarevich--Tate group 
\[ \Sha(E/K) \cong \Sha(E/K)_{\div} \oplus \p{\Sha(E/K)/\Sha(E/K)\div},\]
then from the Cassels--Tate pairing \cite{Cas62} one knows that $\dim_{\bF_p} (\Sha(E/K)/\Sha(E/K)\div)[p]$ is even. As $\dim_{\bF_p} \Sha(E/K)_{\div}[p] = r_p(E/K)$, we obtain

\bthm\lab{ShapCM} If $E$ has complex multiplication defined over $K$, then $\dim_{\bF_p}\Sha(E/K)[p]$ is even. \ethm

\Cref{ShapCM} is a special case of the following weaker version of the Shafarevich--Tate conjecture

\begin{customconj}{$\Sha\mathrm{T}_p(K)$}[\cite{MR10}*{p.~545}]\lab{ShaTp} For every elliptic curve $E/K$, $\dim_{\bF_p} \Sha(E/K)[p]$ is even.  \end{customconj}

For an application of \Cref{ShaTp} to Hilbert's Tenth Problem, see \cite{MR10}*{Theorem 1.2}.

\bigskip
\large{{\it{The strategy of proof of \Cref{G}}}}

\normalsize

It is known (see \S\ref{S18}) that a conjectural formula of T. and V. Dokchitser (\Cref{I}) for the local root number implies \Cref{G} if $p > 2$. In fact, \Cref{F} was proved in \cite{DD08} and \cite{DD11}*{\S5} for $p \ge 3$ by establishing appropriate cases of this formula (recalled in \S\ref{S3}). We introduce it after the following preparations.

\bpp[The setup]\lab{K} Let $K_v$ (to be renamed $K$ from \S\ref{S31} on) be a local field of characteristic $0$ and let $E/K_v$ be an elliptic curve with a ($K_v$-rational) isogeny $E \xra{\phi} E\pr$ of prime degree $p \ge 3$. Let 
\[
\psi \colon \Gal(\ov{K}_v/K_v) \ra \Aut(E[\phi]) \cong \bF_p^\times
\]
be the character giving the Galois action on the kernel $E[\phi]$ of the isogeny and let $K_{v, \psi} = K_v(E[\phi])$ be the fixed field of $\ker \psi$. We denote by $\phi_v\colon E(K_v) \ra E\pr(K_v)$ the map on $K_v$-points induced by $\phi$ and note that the long exact cohomology sequence for $\phi$ tells us that $\#\coker \phi_v$ is finite. Thus it is legitimate to set 
\[ \sigma_{\phi_v} = (-1)^{\ord_p \p{\f{\#\coker \phi_v}{\# \ker \phi_v}}},\]
where $\ord_p a$ is the $p$-adic valuation of $a\in \bQ^\times$. Denoting by $(\psi, -1)_v$ the Artin symbol defined by
\be\lab{J} (\psi, -1)_v = \begin{cases}1, \text{ if } -1 \text{ is a norm in } K_{v, \psi}/K_v, \\ -1, \text{ if not}  \end{cases} \ee
(see \S\ref{U} for another description), we are ready to state \end{pp}

\bconj[$p$-isogeny conjecture \cite{DD11}*{Conjecture 5.3}]\lab{I} $w(E/K_v) = (\psi, -1)_v\cdot \sigma_{\phi_v}$. \econj

\bpp[\Cref{I} implies \Cref{G} for $p \ge 3$] \lab{S18}A well-known global arithmetic duality argument using the Cassels--Poitou--Tate exact sequence (see, for instance, \cite{CFKS10}*{Theorem 2.3}) shows that for $p \ge 3$
\[ (-1)^{\rk_p(E/K)} = \prod_{v \text{ place of } K} \sigma_{\phi_v},\]
 almost all factors on the right hand side being $1$. Hence, \Cref{I}, the product formula for the Artin symbol from global class field theory, and \eqref{B} imply \Cref{G} if $p \ge 3$. \Cref{I} has been settled in many cases, including the case $p = 3$ (see \S\ref{sum} and \Cref{CFKS} for a summary of known cases); to prove \Cref{G} we settle the remaining cases, hence completing the proof of\epp

\bthm[\Cref{DONE}]\lab{II} The $p$-isogeny conjecture is true for $p > 3$. \ethm

\bpp[Progress to date] Using the work of Breuil \cite{Bre00} on classification of finite flat group schemes, Coates, Fukaya, Kato, and Sujatha have proved both a version of \Cref{G} in \cite{CFKS10}*{Theorem 2.1} and a version of \Cref{I} in \cite{CFKS10}*{Theorem 2.7} for abelian varieties of arbitrary dimension. Unfortunately, their methods require a set of hypotheses that exclude some cases of \Cref{I}. For elliptic curves and $p > 3$ their results give\epp

\bthm[\cite{CFKS10}*{Corollary 2.2}]\lab{L} For $p > 3$, the $p$-parity conjecture holds for $E/K$ if $E$ has a ($K$-rational) $p$-isogeny, and if at each place $v$ above $p$ one of the following is true:
\benum[label={(\alph*)}]
\item $E$ has potentially good ordinary reduction at $v$,
\item $E$ has potentially multiplicative reduction at $v$, or
\item $E$ achieves good supersingular reduction after a finite abelian extension of $K_v$.
\eenum \ethm

We make essential use of both the results and the methods of Coates, Fukaya, Kato, and Sujatha in \S\ref{S6} to settle \Cref{I} in the case of a place $v|p$ of additive reduction of Kodaira type III or III$^*$. In all other cases the proof is independent of \Cref{L}.

\bpp[The contents of the paper] In \S\ref{S3} we recall known cases of \Cref{I} and indicate how \Cref{G} was proved by T. and V. Dokchitser if $p \le 3$ (see \S\S\ref{p3}--\ref{p2}). We prove in \S\ref{S5} that the $p$-isogeny conjecture is compatible with making a quadratic twist. The work of \S\ref{S5} is used in \S\ref{S45}, where we settle all the remaining cases of \Cref{I} except those of Kodaira types III or III$^*$. These are taken up in \S\ref{S6} where we make use of the results and methods of Coates, Fukaya, Kato, and Sujatha to finish our proof. In \S\ref{S7} we prove \Cref{CM-Aa,CM-Ba} that concern elliptic curves with complex multiplication.\epp

\bpp[Conventions]\lab{LL} (See also \S\ref{S31} for a notational setup that is valid from \S\ref{S3} on.) Whenever we work with algebraic extensions of a (global or local) field $K$, they are implicitly assumed to lie inside a \emph{fixed} separable closure $\ov{K}$ of $K$. Given a global field $K$ and a place $v$, we implicitly fix an embedding $\ov{K} \hra \ov{K}_v$ and get the corresponding inclusion of a decomposition group $\Gal(\ov{K}_v/K_v) \hra \Gal(\ov{K}/K)$. For a finite flat group scheme (such as $E[\phi]$) over a field $K$ of characteristic 0, we confuse it with its associated Galois representation (such as $E[\phi](\ov{K})$) whenever it is convenient to do so. If $E$ is an elliptic curve over a field $K$, we sometimes write $E/K$ to emphasize the base; if $L/K$ is a field extension, then $E/L$ denotes the corresponding base change $E\times_{\Spec K} \Spec L$. We also make use of the subscript notation to denote base change: for instance, $E[\phi]$ is a $K$-scheme, and $E[\phi]_L$ denotes the base change $E[\phi] \times_{\Spec K} \Spec L$.\epp

\subsection*{Acknowledgements} I thank my advisor Bjorn Poonen for his support and many helpful conversations and suggestions, as well as for reading the manuscript very carefully. I thank Tim Dokchitser for his lectures at the Postech Winter School 2012 in Pohang, South Korea which got me interested in the question answered by this paper, and for pointing out to me \Cref{secret}. Thanks are also due to POSTECH and the organizers of the winter school for an inspiring and hospitable atmosphere. I thank Karl Rubin for telling me that \Cref{CM-Aa} follows from \Cref{G}. I thank Douglas Ulmer for a very helpful conversation about the technique of twisting. I thank Tim Dokchitser, Jessica Fintzen, Jan Nekov{\'a}{\v{r}}, and Bjorn Poonen for comments. I thank the anonymous referee for suggestions and a careful reading of the manuscript.

\section{Known cases of the $p$-isogeny conjecture}\lab{S3}

We have seen in \S\ref{S18} that \Cref{G} follows once we establish the $p$-isogeny conjecture (\Cref{I}). In this section we recall some of the known cases of the latter (\S\S\ref{c1}--\ref{r1}, \S\S\ref{g1}--\ref{UAB}, \S\S\ref{UU}--\ref{p3}), all of which are due to T. and V. Dokchitser \cite{DD08}, \cite{DD11}*{\S5}. In the cases of \S\ref{S} and \S\ref{UA} minor simplifications are provided by \Cref{OMG} and the results of \S\ref{S5}. Since for the rest of the paper we will be working in a local setting, we first change the notation of \S\ref{K} slightly (see also \S\ref{LL} for other conventions). We also recall the classification of local root numbers of elliptic curves in \Cref{RN}. In \S\ref{p2} we discuss the case $p = 2$ which is excluded from \Cref{I}.

\bpp[The setup]\lab{S31} From now on, $K$ denotes a local field of characteristic $0$ (which is assumed to be nonarchimedean from \S\ref{nonarch} on). If $K$ is nonarchimedean, we write $v$ (or $v_K$) for its normalized discrete valuation, $\cO_K$ for the ring of integers, $\fm_K$ for the maximal ideal, $\pi_K$ for a uniformizer, and $\bF_K$ for the residue field. If $L/K$ is a finite extension of nonarchimedean local fields, we write $e_{L/K}$ and $f_{L/K}$ for the ramification index and the degree of the residue field extension, respectively. Let $E/K$ be an elliptic curve with a ($K$-rational) $p$-isogeny $\phi\colon E \ra E\pr$ for a prime $p \ge 3$. The map on $K$-points induced by $\phi$ is denoted by $\phi_K\colon E(K) \ra E\pr(K)$. We write $K_\psi$ for the fixed field of the kernel of the Galois character $\psi\colon \Gal(\ov{K}/K) \ra \Aut(E[\phi]) \cong \bF_p^\times$. We write $(\psi, -1)_K$, instead of $(\psi, -1)_v$, for the Artin symbol \eqref{J}. \Cref{I} becomes
\be\lab{main} 
w(E/K) \overset{?}{=} (\psi, -1)_K \cdot \sigma_{\phi_K}. 
\ee 
\epp

\bthm[\cite{Roh96}*{Theorem 2} and \cite{Kob02}*{Theorem 1.1}, reformulated in \cite{DD10}*{Theorem~3.1}]\label{RN}Let $E$ be an elliptic curve over a local field $K$. Then
\begin{enumerate}[label={(\alph*)}]
\item\lab{RNa} $w(E/K) = -1$, if $K$ is archimedean.
\end{enumerate}
If $K$ is nonarchimedean,
\begin{enumerate}[label={(\alph*)}]\setcounter{enumi}{1}
\item\lab{RNb} $w(E/K) = 1$, if $E$ has good or non-split multiplicative reduction;
\item\lab{RNc} $w(E/K) = -1$, if $E$ has split multiplicative reduction;
\item\lab{RNd} $w(E/K) = \left( \f{-1}{\bF_K}\right)$ if $E$ has additive potentially multiplicative reduction and $\Char \bF_K > 2$ (here $\left( \f{-1}{\bF_K}\right)$ is $1$ if $-1 \in \bF_K^{\times 2}$ and $-1$ otherwise);
\item\lab{RNe} $w(E/K) = (-1)^{\lfloor v(\Delta) \cdot \# \bF_K/12 \rfloor}$, if $E$ has potentially good reduction and $\Char \bF_K > 3$ (here $\Delta$ is a minimal discriminant of $E$).
\end{enumerate}\end{thm}

We now begin the proof of \eqref{main}.

\bpp[The case $K = \bC$]\lab{c1} Since $E[\phi] \subset E(\bC)$ and $\coker \phi_K = 0$, one has $(\psi, -1)_K = 1$ and $\sigma_{\phi_K} = -1$. Since $K$ is archimedean, $w(E/K) = -1$, and \eqref{main} holds.\epp

\bpp[The case $K = \bR$] \lab{r1} Since $\coker \phi_K$ is $2$-torsion, $p\ge 3$, and $\#\ker \phi_K = 1$ or $\#\ker \phi_K = p$, we have $\sigma_{\phi_K} = 1$ or $\sigma_{\phi_K} = -1$, respectively. Accordingly, $K_\psi = \bC$ or $K_\psi = \bR$, so by \eqref{J}, $(\psi, -1)_K = -1$ or $(\psi, -1)_K = 1$. Since $w(E/K) = -1$, \eqref{main} holds in both cases.\epp

Since the archimedean cases of \eqref{main} have been dealt with in \S\S\ref{c1}--\ref{r1}, we assume from now on that $K$ is nonarchimedean.

\bpp\lab{nonarch} If $f \colon A \ra A\pr$ is a homomorphism of abelian groups, we write $\chi(f)$ for $\f{\#\coker f}{\# \ker f}$. Whenever we do so, it is implicitly assumed that the quotient makes sense, i.e., that $\ker f$ and $ \coker f$ are both finite. With this notation, $\sigma_{\phi_K} = (-1)^{\ord_p \chi(\phi_K)}$. We state several elementary properties of $\chi(f)$ that will be used later.\epp

\bprop\lab{M} Given a morphism of short exact sequences
\[ \xymatrix{
0 \ar[r] & A \ar[r]\ar[d]^-{f} & B \ar[r]\ar[d]^-{g} & C\ar[d]^-{h} \ar[r] &0 \\
0 \ar[r] & A\pr \ar[r] & B\pr \ar[r] & C\pr \ar[r] &0, \\
} \]
one has $\chi(g) = \chi(f)\chi(h)$. \eprop

\bpf Snake lemma. \epf

\bprop\lab{O} If $f\colon A \ra A\pr$ and $A, A\pr$ are finite, then $\chi(f) = \f{\# A\pr}{\# A}$. \qed \eprop

\blem\lab{chiP} Let $f\colon K^\times \ra K^\times$ be the $p^{\text{th}}$ power map. Then $\chi(f) = p$ if $\Char \bF_K \neq p$, and $\chi(f) = p^{1 + [K : \bQ_p]}$ if $\Char \bF_K = p$. \elem

\bpf A choice of a uniformizer $\pi_K$ gives $K^\times \cong \bZ \times \cO^\times_K$, so $\chi(f) = p\cdot \chi(\cO_K^\times \xra{f} \cO_K^\times)$ by \Cref{M}. If one uses the filtration of $\cO_K^\times$ by higher units together with the logarithm isomorphism $1 + \fm_K^n \cong \fm_K^n$ for big enough $n$, from \Cref{M,O} one  gets $\chi(\cO_K^\times \xra{f} \cO_K^\times) = \chi(p)$, where $p\colon \fm_K^n \ra \fm_K^n$ is the multiplication by $p$ map. The latter is an isomorphism if $\Char \bF_K \neq p$; if $\Char \bF_K = p$, one has $\chi(p) = p^{[K : \bQ_p]}$.\epf

\bpp[Another description of $\sigma_{\phi_K}$]\lab{R} Let $E_0(K) \subset E(K)$ be the subgroup consisting of points whose reduction lies in the identity component $\wt{\cE}^0$ of the special fiber $\wt{\cE} \ce \cE \times_{\cO_K} \bF_K$ of the N\'{e}ron model $\cE/\cO_K$ of $E/K$. The reduction homomorphism $E(K) \cong \cE(\cO_K) \xra{r} \wt{\cE}(\bF_K)$ is surjective because $\cE/\cO_K$ is smooth and $\cO_K$ is henselian \cite{BLR90}*{\S 2.3 Proposition 5}, so $E(K)/E_0(K) \cong \wt{\cE}(\bF_K)/\wt{\cE}^0(\bF_K)$, which is finite. The order of $E(K)/E_0(K)$ is the \emph{local Tamagawa factor} $c_{E/K}$. The kernel of $r$ is denoted by $E_1(K)$, so $E_0(K)/E_1(K) \cong \wt{\cE}^0(\bF_K)$. The cardinality of the latter is invariant with respect to $K$-rational isogenies: this can be seen, for instance, from the isogeny invariance of the local $L$-factor $L(E/K, s)$ which encodes $\#\wt{\cE}^0(\bF_K)$. Applying \Cref{M,O} to
 \[ \xymatrix{
0 \ar[r] & E_0(K) \ar[r]\ar[d]^-{\phi_0} & E(K) \ar[r]\ar[d]^-{\phi_K} & E(K)/E_0(K)\ar[d] \ar[r] &0 \\
0 \ar[r] & E_0\pr(K) \ar[r] & E\pr(K) \ar[r] & E\pr(K)/E_0\pr(K) \ar[r] &0} \]
and a similar diagram for $E_1(K) \subset E_0(K)$ gives $\chi(\phi_K) = \f{c_{E\pr/K}}{c_{E/K}} \chi(\phi_0) = \f{c_{E\pr/K}}{c_{E/K}} \chi(\phi_1)$. If $p\neq\Char \bF_K$, then $\phi_1$ is an isomorphism \cite{Tat74}*{Corollary 1}, so 
\be\lab{Q} 
\sigma_{\phi_K} = (-1)^{\ord_p \chi(\phi_K)} = (-1)^{\ord_p \f{c_{E\pr/K}}{c_{E/K}}}.
\ee
If $\Char \bF_K = p$, one cannot use this formula. Instead one proceeds as follows: since $E_1(K)$ can be identified with the points of the formal group associated to $E$, one has the canonical exhaustive filtration $E_1(K) \supset E_2(K) \supset \dotsb$ defined by $E_m(K) \ce \ker(\cE(\cO_K) \ra \cE(\cO_K/\fm_K^m))$ (see the proof of \Cref{len2} for another description). The subquotients of the filtration are $E_i(K)/E_{i + 1}(K) \cong \bF_K$ \cite{Tat75}*{\S4}. With this at hand, one applies \Cref{M,O} repeatedly to get $\chi(\phi_1) = \chi(\phi_m)$ for $m \ge 1$. Choose N\'{e}ron differentials $\omega$ on $\cE$ and $\omega\pr$ on $\cE\pr$. Then $\phi^* \omega\pr = \gA \omega$ for some $\gA\in \cO_K$. For $m$ large enough, formal logarithm furnishes isomorphisms $ E_m(K) \isomto \fm_K^m$, $E_m\pr(K) \isomto \fm_K^m$, under which, by \cite{Rub99}*{Proposition 3.14}, $\phi_m$ corresponds to multiplication by $\gA$. The analogue of \eqref{Q} is therefore
\be\lab{QQ} 
\sigma_{\phi_K} = (-1)^{v(\gA)f_{K/\bQ_p} + \ord_p \f{c_{E\pr/K}}{c_{E/K}}}.
\ee 
Formula \eqref{QQ} is a special case of a more general formula \cite{Sch96}*{Lemma 3.8} valid for abelian varieties of any dimension. 

If $E$ has potentially good reduction, the local Tamagawa factors are at most $4$ \cite{Tat74}*{Addendum to Theorem 3}. If in addition $p > 3$, then \eqref{Q} for the case $\Char \bF_K \neq p$ becomes
\[
\sigma_{\phi_K} = 1,
\]
and \eqref{QQ} for the case $\Char \bF_K = p$ becomes 
\be\lab{QQQ} 
\sigma_{\phi_K} = (-1)^{v(\gA)f_{K/\bQ_p}}.
\ee
\epp

\bpp[The case $\Char \bF_K \neq p$, and $E$ has good reduction]\lab{g1} The local Tamagawa factors are $1$, so \eqref{Q} gives $\sigma_{\phi_K} = 1$. Also, $K_\psi \subset K(E[p])$, and the latter is an unramified extension of $K$ by the criterion of N\'{e}ron--Ogg--Shafarevich. Hence, so is $K_\psi/K$, and from \eqref{J} we get $(\psi, -1)_K = 1$. By \Cref{RN}~\ref{RNb}, $w(E/K) = 1$, and \eqref{main} holds.\epp

\bpp[The case $\Char \bF_K \neq p$, and $E$ has potentially multiplicative reduction]\lab{S} If the reduction is not split multiplicative, it becomes so after a unique quadratic extension $L/K$ which is unramified if and only if $E$ has multiplicative reduction \cite{Tat74}*{pp.~190--191}. The quadratic twist of $E$ by $L/K$ has split multiplicative reduction (this can be seen, for instance, from \Cref{LLL}). By \Cref{twist}, we are therefore reduced to the case of split multiplicative reduction (see \Cref{ind}).

Using Tate's theory of rigid analytic uniformization (loc.~cit.) one has $E \cong \bG_m/q^\bZ$ for some $q\in K^\times$ with $v(q) > 0$ (and similarly for $E\pr$). Using rigid-analytic GAGA one sees that $\phi\colon E\ra E\pr$ can be written as $\bG_m/q^\bZ \ra \bG_m/(q^p)^\bZ$ induced by the $p^{\text{th}}$ power map, or $\bG_m/(q^{\prime p})^\bZ \ra \bG_m/(q^{\prime})^\bZ$ induced by the identity. By \Cref{M} and \Cref{chiP}, in the first case $\chi(\phi_K) = p$; in the second $\chi(\phi_K) = 1/p$. In both cases $\sigma_{\phi_K} = -1$. In the first case $K_\psi = K(\zeta_p)$; in the second $K_\psi = K$. In both cases $K_\psi/K$ is unramified, so $(\psi, -1)_K = 1$. \Cref{RN}~\ref{RNc} gives $w(E/K) = -1$, and \eqref{main} holds. \epp

\bpp[The case $\Char \bF_K \neq p$, $p > 3$, and $E$ has additive potentially good reduction]\lab{UAB} This is \cite{DD08}*{Lemma 9}. \epp

\bpp[Another description of the Artin symbol]\lab{U} For a continuous character 
\[ \theta\colon \Gal(\ov{K}/K) \ra \bQ/\bZ\]
and $a \in K^\times$, let $(\theta, a)_K \in \bQ/\bZ$ be the corresponding symbol (cf.~\cite{Ser79}*{Chapter XIV, \S 1}). It can also be defined by setting $(\theta, a)_K \ce \theta(\Art_K(a))$ where $\Art_K \colon K^\times \ra \Gal(\ov{K}/K)^{\ab}$ is the local Artin homomorphism. Therefore, the pairing $(\theta, a)_K$ is bilinear, and if $K_\theta$ is the fixed field of $\ker \theta$, then $(\theta, a)_K$ vanishes if and only if $a$ is a norm in $K_\theta/K$. Since we are only interested in the case $a = -1$, we think of $(\theta, -1)_K$ as taking values in $\{ \pm 1\}$. 

Fix an injection $\bF_p^\times \hra \bQ/\bZ$. One possible choice for $\theta$ is $\psi$ (\S\ref{S31}). In this case $K_\psi = K(E[\phi])$, and we recover the Artin symbol $(\psi, -1)_K$ (cf.~\eqref{J}). Another possible choice is the cyclotomic character $\omega: \Gal(\ov{K}/K) \ra \bF_p^\times$ described as follows: for $\zeta_p$ a primitive $p^{\text{th}}$ root of unity $s(\zeta_p) = \zeta_p^{\omega(s)}$ (this is independent of $\zeta_p$). In this case $K_\omega = K(\zeta_p)$.\epp

\begin{lemma}\lab{Artin} For a finite extension $L/K$ and every $\theta$ as above, one has $(\theta, -1)_L = (\theta, -1)_K^{[L:K]}$.\end{lemma}
\begin{proof} This follows from the bottom diagram of \cite{Ser67}*{\S2.4} applied to $-1$. \epf

\bcor\lab{OMG} Suppose that $K$ is a finite extension of $\bQ_p$. Then $(\omega, -1)_K = (-1)^{[K:\bQ_p]}$. \ecor

\bpf Since $(\omega, -1)_{\bQ_p} = -1$ (see, e.g., \cite{Ser67}*{\S3.1 Theorem 2 (2)}), one applies \Cref{Artin}. \epf

To handle the cases of \eqref{main} when $\Char \bF_K = p$ we need to be able to compute symbols $(\theta, -1)_K$ for continuous characters $\theta\colon \Gal(\ov{K}/K) \ra \bF_p^\times$, in which case $K_\theta$ is a cyclic extension of $K$ of degree dividing $p - 1$. These symbols have been worked out by T. and V. Dokchitser:

\blem[\cite{DD08}*{Lemma 12}]\label{JNT} Let $K$ be a finite extension of $\bQ_p$ with $p$ odd, and fix a character $\theta\colon \Gal(\ov{K}/K) \ra \bF_p^\times$. Then \[ (\theta, -1)_K = (-1)^{f_{K/\bQ_p}(p - 1)/e_{K_\theta/K}};\] in other words, $(\theta, -1)_K = 1$ if and only if either $f_{K/\bQ_p}$ or $\f{p - 1}{e_{K_\theta/K}}$ is even.
 \elem

\bpp[The case $\Char \bF_K = p > 3$, $E$ has potentially good reduction, and $f_{K/\bQ_p}$ is even]\lab{UU} We show that all the terms in \eqref{main} are $1$. \Cref{JNT} gives $(\psi, -1)_K = 1$, whereas \eqref{QQQ} gives $\sigma_{\phi_K} = 1$. To compute $w(E/K)$ we use \Cref{RN}~\ref{RNe}: $\#\bF_K$ is a square, so $\# \bF_K \equiv 1 \bmod 24$, whereas $v(\Delta) < 12$ because the reduction is potentially good. Hence, $w(E/K) = (-1)^{\lfloor v(\Delta)/12 \rfloor} = 1$.\epp

\bpp[The case $\Char \bF_K = p$, and $E$ has good reduction]\lab{U0} Because of \S\ref{UU}, we assume that $f_{K/\bQ_p}$ is odd. By \Cref{RN}~\ref{RNb}, $w(E/K) = 1$, so by \Cref{JNT} and \eqref{QQQ}, to check \eqref{main} one needs to argue that
\[\lab{ten} (-1)^{(p - 1)/e_{K_\psi/K}} = (-1)^{v(\gA)}, \]
where $\gA$ defined in \S\ref{R} is also the coefficient of $T$ in the power series $f(T)$ giving the action of $\phi$ on formal groups. This is done in \cite{DD08}*{\S6} by a careful analysis of $f(T)$.\epp

\bpp[The case $\Char \bF_K = p$, and $E$ has potentially multiplicative reduction]\lab{UA} The argument is the same as in \S\ref{S}, and the second case there requires no modification. In the first case, by \Cref{M} and \Cref{chiP}, $\chi(\phi_K) = p^{1 + [K : \bQ_p]}$, so $\sigma_{\phi_K} = (-1)^{1 + [K : \bQ_p]}$. We have $K_\psi = K(\zeta_p)$, so by \eqref{J} and \Cref{OMG}, $(\psi, -1)_K = (\omega, -1)_K = (-1)^{[K : \bQ_p]}$. By \Cref{RN}~\ref{RNc}, $w(E/K) = -1$, so \eqref{main} holds.
\epp

\bpp[The case $p = 3$]\lab{p3} The argument of \cite{DD08}*{Lemma 9} used in \S\ref{UAB} faces complications if $p = 3$. Calculations are still manageable if $\Char \bF_K \neq p$ (see \cite{DD08}*{Lemma 10}), but get out of hand if $\Char \bF_K = p$. To treat these cases, and hence establish \eqref{main} if $p = 3$, T. and V. Dokchitser have used a global-to-local deformation argument with \cite{Nek09}*{Theorem 1} as an input, see \cite{DD11}*{Theorem~5.7}. \epp

\bpp[The case $p = 2$]\lab{p2} The formula \eqref{main} does not hold if $p = 2$. Indeed, the kernel of a $2$-isogeny is always rational, so $K_\psi = K$, giving $(\psi, -1)_K = 1$. If $K = \bR$ and $\phi_K$ is not surjective, $\#\coker \phi_K = 2$, and $\sigma_{\phi_K} = 1$. This violates \eqref{main}, because by \Cref{RN}~\ref{RNa} $w(E/\bR) = -1$. 

In case $p = 2$, \Cref{G} was proved by T. and V. Dokchitser by finding an analogue \cite{DD11}*{Conjecture 5.3} of \eqref{main}, proving it in most cases by direct computations in \cite{DD08}*{\S7}, and using a global-to-local deformation argument to handle the remaining ones in \cite{DD11}*{Theorem 5.7}.\epp

\bpp[Summary of the known cases]\lab{sum} By \S\S\ref{c1}--\ref{r1}, \S\S\ref{g1}--\ref{UAB}, and \S\S\ref{UU}--\ref{p3}, the $p$-isogeny conjecture \eqref{main} is true in all the cases, except possibly when $\Char \bF_K = p > 3$, $f_{K/\bQ_p}$ is odd, and the reduction is additive potentially good. After some preparations in \S\ref{S5} we take up the remaining cases in \S\ref{S45}.\epp




\section{Compatibility with making a quadratic twist}\lab{S5}

To verify that the $p$-isogeny conjecture \eqref{main} is compatible with making a quadratic twist (\Cref{twist}), we investigate how its individual terms change under this operation (\Cref{X,Y,DD}). The technique of twisting is standard but to fix ideas we begin by recalling the way in which we prefer to think about it. The setup is that of \S\ref{S31}.

\bpp[Twisting by a quadratic Galois character]\lab{S51} Suppose that $L/K$ is a quadratic extension and let $\tau\colon \Gal(\ov{K}/K) \ra \{ \pm 1\}$ be the corresponding nontrivial character. Since $\{ \pm 1\} \le \Aut_{\ov{K}}(E)$ we can think of $\tau\colon \Gal(\ov{K}/K) \ra \Aut_{\ov{K}}(E)$ as a (crossed) homomorphism and therefore identify $\tau$ with the corresponding element of $H^1(K, \Aut_{\ov{K}}(E))$. The elements of the latter pointed set classify \emph{twists} of $E$ (cf.~\cite{Ser02}*{ I.\S5.3} or \cite{BS64}*{Proposition 2.6}), i.e., elliptic curves over $K$ that are $\ov{K}$-isomorphic to $E$. In particular, $\tau$ gives rise to the \emph{twist} $\wt{E}/K$ of $E/K$ by $L/K$. \epp

For a $K$-scheme $X$, the $X$-valued points of $\wt{E}$ are the $X\times_K \ov{K}$-valued points of $E_{\ov{K}}\ce E \times_K \ov{K}$ that are invariant under the twisted by $\tau$ Galois action, i.e., those $P \in E(X\times_K \ov{K})$ with $P = \tau(s)\cdot  \leftexp{s}P$ for all $s\in \Gal(\ov{K}/K)$. These are $P \in E(X\times_K L)$ on which the nontrivial element $t \in \Gal(L/K)$ acts as $\leftexp{t}P = -P$. This gives a description of the functor of points of $\wt{E}/K$; we also see that the isomorphism $E_{\ov{K}} \cong \wt{E}_{\ov{K}}$ is $\Gal(\ov{K}/L)$-equivariant, so $E$ and $\wt{E}$ are $L$-isomorphic. Twisting being functorial in $E$ (loc.~cit.), $\wt{E}$ possesses a $p$-isogeny $\wt{\phi}\colon \wt{E} \ra \wt{E\pr}$ defined over $K$.

\brem If $L = K(\sqrt{d})$ and one wishes to think in terms of Weierstrass equations 
\be\lab{eqE} 
E\colon y^2 = x^3 + Ax + B,
\ee
then the quadratic twist described above is the usual $\wt{E}\colon dy^2 = x^3 + Ax + B$. Indeed, since $[-1]_E$ in these coordinates is $(x, y) \mapsto (x, -y)$, multiplying the $y$-coordinate by $\sqrt{d}$ has the desired effect as far as the Galois action is concerned. By scaling the variables one can bring the equation for $\wt{E}$ to 
\be\lab{eqEE} y^2 = x^3 + Ad^2x + Bd^3.\ee
The discriminants $\Delta$ and $\wt{\Delta}$ of \eqref{eqE} and \eqref{eqEE} are related by
\be\lab{DDD} 
\wt{\Delta} = d^6 \Delta. 
\ee 
\erem

\bpp[Implications for the $p$-isogeny] Consider the restriction of scalars $N_{L/K} E$ of $E/L$ back to $K$. By definition, $(N_{L/K} E)(X) = E(X\times_K L)$ functorially in $X$ and $E$. We have seen that $\wt{E}(X)$ identifies with the $-1$-eigenspace of $(N_{L/K} E)(X)$ for the action of $t$; similarly $E(X)$ identifies with the $+1$-eigenspace. These identifications being functorial, one gets a $K$-homomorphism of abelian varieties 
\[ f_E: E\times \wt{E} \ra N_{L/K} E.\]
Since the intersection of the eigenspaces consists of $2$-torsion points, so does the kernel of $f_E$. We conclude that $f_E$ is an isogeny and that
\[ \p{E(X) \tensor \bZ\left[\f{1}{2}\right]} \oplus \p{\wt{E}(X) \tensor \bZ\left[\f{1}{2}\right]} \isomto (N_{L/K} E)(X) \tensor \bZ\left[\f{1}{2}\right].\]
The latter being functorial in $X$ and $E$, taking $X = \Spec K$ we get the commutative diagram
\be\lab{W}\begin{split} \xymatrix{
 (E(K)\tensor \bZ[\f{1}{2}]) \oplus (\wt{E}(K)\tensor \bZ[\f{1}{2}])\ar[d]_-{(\phi_K \tensor \bZ[\f{1}{2}]) \oplus (\wt{\phi}_K\tensor \bZ[\f{1}{2}])
} \ar[r]^-{\sim} &E(L)\tensor \bZ[\f{1}{2}] \ar[d]^-{\phi_L \tensor \bZ[\f{1}{2}]} \\
 (E\pr(K)\tensor \bZ[\f{1}{2}]) \oplus (\wt{E\pr}(K)\tensor \bZ[\f{1}{2}]) \ar[r]^-{\sim} &E\pr(L)\tensor \bZ[\f{1}{2}]. 
}\end{split}\ee
Since $- \tensor \bZ[\f{1}{2}]$ is exact and does not affect the $p$-primary parts ($p > 2$), \eqref{W} gives \epp

\begin{prop}\lab{X}$\sigma_{\phi_K}\sigma_{\wt{\phi}_K} = \sigma_{\phi_L}.$\end{prop}

\bpp[The twist of $\psi$] Describing the character $\wt{\psi}\colon \Gal(\ov{K}/K) \ra \bF_p^\times$ that gives the Galois action on $\wt{E}[\wt{\phi}]$ is easy: $\wt{E}[\wt{\phi}](\ov{K})$ identifies with $E[\phi](\ov{K})$ with the Galois action twisted by $\tau$, so $\wt{\psi} = \psi\tau$.\epp

\begin{prop}\lab{Y}$(\psi, -1)_L = 1 = (\psi, -1)_K(\wt{\psi}, -1)_K(\tau, -1)_K$.\end{prop}
\pf To deal with the left hand side one applies \Cref{Artin}. The right hand side is taken care of by bilinearity (\S\ref{U}):
\[ (\psi, -1)_K(\wt{\psi}, -1)_K(\tau, -1)_K = (\psi, -1)_K^2 (\tau, -1)_K^2 = (\psi\tau, 1)_K = 1.\qedhere \]
\epf

\bpp[Implications for the $l$-adic Tate module] The isogeny $f_E$ induces an isomorphism of $l$-adic Tate modules ($l\neq 2$) 
\[ 
V_l(E/K) \oplus V_l(\wt{E}/K) \cong V_l((E\times \wt{E})/K) \cong V_l((N_{L/K}E)/K) \cong \Ind_L^K V_l(E/L).
\]
The last isomorphism is a general property of the restriction of scalars for abelian varieties (use \cite{Mil72}*{Proposition 6 (b)} together with the projection formula $(\Ind_L^K \b{1}_L) \tensor V_l(E/K) \cong \Ind_L^K V_l(E/L)$). Of course, this gives an isomorphism of the corresponding Weil--Deligne representations \cite{Roh94}*{\S3--\S4 and \S13} (take $l \neq 2, \Char \bF_K$)
\be\lab{BB} \sigma\pr_{E/K} \oplus \sigma\pr_{\wt{E}/K} \cong \Ind_L^K \sigma\pr_{E/L}. \ee\epp

\begin{prop}\lab{LLL} The local $L$-factors are related by \[ L(E/K, s)L(\wt{E}/K, s) = L(E/L, s).\] \eprop
\bpf This follows immediately from \eqref{BB} and the multiplicativity and inductivity of the $L$-factor of a Weil--Deligne representation \cite{Roh94}*{\S8 and \S17}. \epf

\bpp[Properties of local root numbers]\lab{EE} Let $\eta\colon K \ra \bC^\times$ be a nontrivial (continuous) additive character, and let $dx$ be a Haar measure on $(K, +)$. Let $\sigma\pr$ be a Weil--Deligne representation, and let $\eps(\sigma\pr, \eta, dx)$ be its $\eps$-factor (cf.~\cite{Roh94}). The \emph{root number} of $\sigma\pr$ is
\[ w(\sigma\pr, \eta) \ce \f{\eps(\sigma\pr, \eta, dx)}{\abs{\eps(\sigma\pr, \eta, dx)}}.\]
Standard properties of $\eps$-factors show (op.~cit.)~that $w(\sigma\pr, \eta)$ is independent of the choice of $dx$ and is even independent of $\eta$ if $\sigma\pr$ is essentially symplectic, in which case we write $w(\sigma\pr)$. Due to Weil pairing, this is the case for $\sigma\pr_{E/K}$ associated to an elliptic curve $E$ (op.~cit.);~by definition $w(E/K) = w(\sigma\pr_{E/K})$. For use in the proof of \Cref{DD} we record some basic properties of the root number, all of which follow from analogous properties of the $\eps$-factor (op.~cit.):
\begin{enumerate}[label=(\alph*)]
\item\lab{EEa} Additivity: $w(\sigma\pr_1 \oplus \sigma_2\pr, \eta) = w(\sigma_1\pr, \eta)w(\sigma_2\pr, \eta)$;
\item\lab{EEb} Inductivity in degree zero: if $L/K$ is a finite extension and $\sigma\pr$ is a virtual representation of degree $0$ of the Weil--Deligne group $\cW\pr(\ov{K}/L)$ of $L$, then $w(\Ind_L^K \sigma\pr, \eta ) = w(\sigma\pr, \eta\circ \tr_{L/K})$;

\item\lab{EEc} 
Determinant formula: $w(\sigma\pr, \eta)w(\sigma^{\prime *}, \eta) = (\det \sigma)(-1)$, where $\sigma^{\prime *}$ is the contragredient of $\sigma\pr$, $\sigma$ is the underlying representation of the Weil group $\cW(\ov{K}/K)$, and, with the local Artin homomorphism $\Art_K$ from \S\ref{U}, $(\det \sigma)(-1) = \det \sigma \circ \Art_K(-1)$.
\end{enumerate}
In particular, applying \ref{EEc}~to a self-contragredient character such as $\b{1}_K$ or $\tau$ we get 
\be\lab{EEE} w(\b{1}_K, \eta)^2 = 1\quad\text{ and }\quad w(\tau, \eta)^2 = \tau(-1) = (\tau, -1)_K.\ee\epp

\blem\lab{Eind} $w(\Ind_L^K \b{\upshape{1}}_L, \eta)^2 = (\tau, -1)_K.$ \elem

\bpf Using decomposition $\Ind_L^K \b{1}_L \cong \b{1}_K \oplus \tau$, we compute
\[ w(\Ind_L^K \b{1}_L, \eta)^2 \overset{\S\ref{EE}~\ref{EEa}}{=} w(\b{1}_K, \eta)^2 w(\tau, \eta)^2 \overset{\eqref{EEE}}{=} (\tau, -1)_K.\qedhere  \] 
\epf

\begin{prop}\lab{DD}The local root numbers are related by \[ w(E/L) = w(E/K)w(\wt{E}/K)(\tau, -1)_K.\]\end{prop}
\pf From \eqref{BB} we get
\[\ba w(E/K)w(\wt{E}/K) &= w(\Ind_L^K \sigma\pr_{E/L}, \eta) \overset{\S\ref{EE}~\ref{EEb}}{=} \f{w(E/L)}{w(\b{1}_L \oplus \b{1}_L, \eta \circ \tr_{L/K})}w(\Ind_L^K (\b{1}_L \oplus\b{1}_L), \eta)\\ & \overset{\eqref{EEE}}{=} w(E/L)w(\Ind_L^K \b{1}_L, \eta)^2 \overset{\ref{Eind}}{=} w(E/L)(\tau, -1)_K. \ea\]
Since $(\tau, -1)_K \in \{ \pm 1\}$, we can carry it to the other side, and the conclusion follows.
\epf

\begin{thm}\lab{twist}Fix a prime $p \ge 3$. The $p$-isogeny conjecture \eqref{main} is compatible with quadratic twists: in the setup of \S\ref{S31} and \S\ref{S51},
\[ w(E/L)\cdot ((\psi, -1)_K\cdot \sigma_{\phi_{K}})\cdot ((\wt{\psi}, -1)_K \cdot \sigma_{\wt{\phi}_{K}}) = ((\psi, -1)_L\cdot \sigma_{\phi_L}) \cdot  w(E/K)\cdot w(\wt{E}/K). \]
In particular, if the $p$-isogeny conjecture holds for two of $E/K, E/L, \wt{E}/K$, it holds for the third one. \end{thm}

\bpf Combine \Cref{X,Y,DD}. \epf

\brem\lab{ind} Note that \Cref{twist} holds regardless of $\Char \bF_K$, and in its proof we have not used the case-by-case analysis of the $p$-isogeny conjecture from \S\ref{S3}. In particular, it was legitimate to use it in \S\ref{S} and \S\ref{UA}. \erem

\section{The remaining cases of the $p$-isogeny conjecture}\lab{S45}

We have seen in \S\ref{sum} that the $p$-isogeny conjecture holds in most cases, including all cases when $\Char \bF_K \neq p$. In this section we prove it in all of the remaining cases except for Kodaira types III or III$^*$, which are treated in \S\ref{S6}.

\bpp[The restricted setup]\lab{S61} Since the $p$-isogeny conjecture is known in other cases (see \S\ref{sum}), for the rest of the paper we make the following assumptions in addition to those of \S\ref{S31}: the degree of the isogeny is equal to the residue characteristic $\Char \bF_K = p > 3$, the reduction of $E$ is additive potentially good, and the degree of the residue field extension $f_{K/\bQ_p}$ is odd. Let $\Delta$ be a minimal discriminant of $E$. Define $\Delta\pr$ similarly for $E\pr$. Since $\Char \bF_K > 3$ and the reduction is potentially good, $v(\Delta) < 12$ \cite{Tat75}*{p.~46}. In fact, $v(\Delta) = 2, 3, 4, 6, 8, 9$, or $10$, corresponding to Kodaira types II, III, IV, I$_0^*$, IV$^*$, III$^*$, or II$^*$, respectively (loc.~cit.).\epp

\begin{lemma}\lab{HH}Suppose that $L/K$ is a finite extension of ramification index $e = e_{L/K}$, and that the degree $f_{L/K}$ of the residue field extension is odd. Write $ev(\Delta) = 12b + a$ with $0 \le a < 12$, so that $b =  \lfloor e v(\Delta)/12 \rfloor$. Define $a\pr, b\pr$ analogously using $\Delta\pr$. Then 
\begin{enumerate}[label={(\alph*)}]
\item\lab{HHa} $a$ and $a\pr$ are the $L$-valuations of minimal discriminants of $E/L$ and $E\pr/L$, respectively, and 
\item\lab{HHb} $\sigma_{\phi_L} = \sigma_{\phi_K}^e \cdot (-1)^{b + b\pr}.$
\eenum \end{lemma}

\bpf Choose minimal equations for $E/K$ and $E\pr/K$ to get associated minimal discriminants $\Delta$ and $\Delta\pr$, and N\'{e}ron minimal differentials $\omega$ and $\omega\pr$. When we pass from $K$ to $L$ those equations might not be minimal anymore: one may need to make changes of coordinates to arrive at minimal equations over $L$. When making those changes of coordinates one will have some $u, u\pr \in L$ for which $\Delta, \Delta\pr$ will get multiplied by $u^{-12}, (u\pr)^{-12}$, respectively, and $\omega, \omega\pr$ will get multiplied by $u, u\pr$, respectively \cite{Del75}*{(1.2) et (1.8)}. Since the reduction will stay potentially good, the $L$-valuations of new minimal discriminants will be $< 12$ and will therefore equal $a$ and $a\pr$, respectively, giving \ref{HHa}. Also, $v_L(u) = b$, $v_L(u\pr) = b\pr$, and \ref{HHb} follows from \eqref{QQQ} because we are assuming that $f_{K/\bQ_p}$ and $f_{L/\bQ_p}$ are odd.\epf

\begin{rem}\lab{MM}The set $\{ v(\Delta), v(\Delta\pr) \}$ is a subset of one of the following: $\{2, 10\}$, $\{3, 9\}$, $\{4, 8\}$, or $\{6\}$. This is because $E$ acquires good reduction over an extension $L/K$ if and only if $E\pr$ does, whereas \Cref{HH}~\ref{HHa} tells us that the minimal ramification index of an extension over which $E$ acquires good reduction is $\f{12}{\gcd(v(\Delta), 12)}$. Hence, $\gcd(v(\Delta), 12) = \gcd(v(\Delta\pr), 12)$.\end{rem}

\blem\lab{TTT} Let $L/K$ be a ramified quadratic extension, and let $\wt{E}$ be the corresponding twist of $E$. Let $\wt{\Delta}$ be a minimal discriminant of $\wt{E}$. Then $v(\wt{\Delta}) \equiv v(\Delta) + 6 \bmod 12$. \elem

\bpf Since $\Char \bF_K > 3$ and the reduction of $\wt{E}$ is potentially good, $v(\wt{\Delta}) < 12$ \cite{Tat75}*{p.~46}. But $L = K(\sqrt{\pi_K})$ for some uniformizer $\pi_K\in \cO_K$, so the conclusion follows from \eqref{DDD}.\epf

\bprop\lab{six} Under the assumptions of \S\ref{S61}, the $p$-isogeny conjecture \eqref{main} is true if $v(\Delta) = 6$. \eprop

\bpf \Cref{HH}~\ref{HHa} shows that $E$ acquires good reduction after a quadratic ramified extension $L/K$. The corresponding quadratic twist has good reduction by \Cref{TTT}. The conclusion then follows from \Cref{twist} and \S\ref{U0}.  \epf

The following relation between the discriminants of elliptic curves related by a $p$-isogeny has been communicated to me by Tim Dokchitser:

\bthm[\cite{Coa91}*{Appendix, Theorem 8}, \cite{DD12}*{Theorem 1.1}]\lab{secret} Let $E$ and $E\pr$ be elliptic curves over a field $K$ of characteristic $0$. Suppose that $\phi\colon E \ra E\pr$ is a $p$-isogeny with $p > 3$. Let $\Delta$ and $\Delta\pr$ be discriminants of some Weierstrass equations for $E$ and $E\pr$, respectively. Then $\Delta\pr/\Delta^p \in (K^\times)^{12}$ regardless of the Weierstrass equations chosen. \ethm

\blem\lab{discs} Under the assumptions of \S\ref{S61}, suppose that $v(\Delta) = 4$ or $v(\Delta) = 8$. Then $v(\Delta) = v(\Delta\pr)$ if and only if $\#\bF_K \equiv 1 \bmod 6$. \elem

\bpf 
\Cref{secret} tells us that $v(\Delta\pr) - p\cdot v(\Delta) \equiv 0 \bmod 12$. Hence, if $v(\Delta) = v(\Delta\pr)$, then $p - 1$ is divisible by $3$, so $\#\bF_K = p^{f_{K/\bQ_p}} \equiv p \equiv 1 \bmod 6$. Conversely, if $\#\bF_K \equiv 1 \bmod 6$, then $p \equiv 1 \bmod 6$, and hence $v(\Delta\pr) \equiv v(\Delta)\bmod 6$, so by \Cref{MM} we must have $v(\Delta) = v(\Delta\pr)$.
\epf

\bprop\lab{four} Under the assumptions of \S\ref{S61}, the $p$-isogeny conjecture \eqref{main} is true if $v(\Delta) = 4$ or $v(\Delta) = 8$. \eprop

\pf Choose a cubic totally ramified extension $L/K$. By \Cref{HH}~\ref{HHa}, $E/L$ has good reduction, so by \S\ref{U0} it satisfies the $p$-isogeny conjecture. We check how the terms change when passing from $K$ to $L$:
\begin{enumerate}[label={(\alph*)}]
\item $(\psi, -1)_L = (\psi, -1)_K$ by \Cref{Artin}.
\item The valuation $v(\Delta)$ is $4$ or $8$, and \Cref{HH} implies that $b$ is $1$ or $2$ accordingly. By \Cref{MM}, $v(\Delta\pr)$ also is $4$ or $8$, so $b\pr$ is $1$ or $2$ accordingly. Thus, by \Cref{HH}~\ref{HHb}, $\sigma_{\phi_L} \neq \sigma_{\phi_K}$ if and only if $v(\Delta) \neq v(\Delta\pr)$. By \Cref{discs}, this is the case if and only if $\# \bF_K \equiv 5 \bmod 6$.

\item By \Cref{RN}~\ref{RNe}, $w(E/L) = 1$, whereas $w(E/K) = 1$ if $\# \bF_K \equiv 1 \bmod 6$, and $w(E/K) = -1$ if $\# \bF_K \equiv 5 \bmod 6$. Therefore, $w(E/L) \neq w(E/K)$ if and only if $\#\bF_K \equiv 5 \bmod 6$.
\end{enumerate}  
We conclude that when passing from $K$ to $L$ both sides of \eqref{main} change sign if and only if $\#\bF_K \equiv 5 \bmod 6$. Since \eqref{main} holds for $E/L$ by \S\ref{U0}, it must hold for $E/K$ as well.\epf

\begin{prop}\lab{two} Under the assumptions of \S\ref{S61}, the $p$-isogeny conjecture \eqref{main} is true if $v(\Delta) = 2$ or $v(\Delta) = 10$.\end{prop}

\pf Choose a quadratic ramified extension $L/K$. By \Cref{HH}~\ref{HHa}, the valuation of a minimal discriminant of $E/L$ is $4$ or $8$. By \Cref{TTT}, the valuation of a minimal discriminant of the twist $\wt{E}$ is $8$ or $4$. In particular, the $p$-isogeny conjecture holds for both $E/L$ and $\wt{E}/K$ by \Cref{four}. By \Cref{twist}, it must hold for $E$ as well.\epf

\brem Another way to prove \Cref{two} is to choose a ramified cubic extension $L/K$ and check that none of the terms in \eqref{main} change when passing from $K$ to $L$. The argument is similar to that of \Cref{7A}. \erem

\section{The case of Kodaira type III or III$^*$}\lab{S6}

As pointed out in \S\ref{S61}, this is the case when $v(\Delta) = 3$ or $v(\Delta) = 9$. To study it we are going to use the work of Coates, Fukaya, Kato, and Sujatha \cite{CFKS10} that settles the $p$-isogeny conjecture in many cases. We begin with a lemma that will be useful later in imposing additional assumptions in \Cref{7B}.

\blem\lab{7A} Under the assumptions of \S\ref{S61}, suppose that $v(\Delta) = 3$ or $v(\Delta) =  9$ and let  $L/K$ be an extension of odd degree. The $p$-isogeny conjecture \eqref{main} holds for $E/K$ if and only if it holds for $E/L$.\elem

\bpf In fact, none of the terms in \eqref{main} change when passing from $K$ to $L$:
\begin{enumerate}[label={(\alph*)}]
\item $(\psi, -1)_L = (\psi, -1)_K$ by \Cref{Artin}.
\item\lab{7Ab} $\sigma_{\phi_L} = \sigma_{\phi_K}$ because in \Cref{HH}~\ref{HHb} one has $b \equiv b\pr \bmod 2$. Indeed, by \Cref{MM}, $v(\Delta\pr) \in \{3, 9\}$, so one only needs to check that $\lfloor 3 \cdot e_{L/K}/12 \rfloor \equiv \lfloor 9 \cdot e_{L/K}/12 \rfloor \bmod 2$, or equivalently that $\lfloor  e_{L/K}/4 \rfloor \equiv \lfloor 3 \cdot e_{L/K}/4 \rfloor \bmod 2$. This is confirmed after a short check of possibilities $e_{L/K} \in \{1, 3, 5, 7\} \bmod 8$.
\item The valuation $v_L(\Delta_L)$ of a minimal discriminant of $E/L$ is in $ \{ 3, 9\}$ by \Cref{HH}~\ref{HHa}. Also, \Cref{RN}~\ref{RNe} yields $w(E/L) = (-1)^{\lfloor v_L(\Delta_L) \cdot \#\bF_L/12 \rfloor}$. The latter is $(-1)^{\lfloor v_L(\Delta_L) \cdot \#\bF_K/12 \rfloor}$, because $\#\bF_L \equiv \#\bF_K \bmod 8$. Since $w(E/K) = (-1)^{\lfloor v(\Delta) \cdot \#\bF_K/12 \rfloor }$, to check that $w(E/L) = w(E/K)$ one needs to check that $\lfloor 3 \cdot \# \bF_K /12 \rfloor \equiv \lfloor 9 \cdot \# \bF_K /12 \rfloor \bmod 2$, which is the same computation as in \ref{7Ab}. \qedhere
\end{enumerate} \epf

\bthm[\cite{CFKS10}*{Theorem 2.7 and Proposition 2.8 (3)}]\lab{CFKS} Assume the setup of \S\ref{S61}. The $p$-isogeny conjecture \eqref{main} holds if either $E$ has potentially good ordinary reduction, or $E$ achieves good supersingular reduction over a finite abelian extension of $K$. \ethm

\bpp[Consequences for the case at hand]\lab{S70} Assume the setup of \S\ref{S61} and suppose that $v(\Delta) = 3$ or $v(\Delta) = 9$. Let $F/K$ be a totally ramified extension of degree $4$, so $E/F$ has good reduction by \Cref{HH}~\ref{HHa}. If $K$ contains a primitive $4^{\text{th}}$ root of unity, i.e., if $-1 \in \bF_K^{\times 2}$, the extension $F/K$ is abelian and we can apply \Cref{CFKS} to deduce \eqref{main}. If $-1 \not\in \bF_K^{\times 2}$, then $\# \bF_K \equiv 3 \bmod 4$, and because of \Cref{CFKS} we can assume in addition that $E/K$ is potentially supersingular.  \epp

\blem\lab{7B} Under the assumptions of \S\ref{S61}, suppose that $v(\Delta) = 3$ or $v(\Delta) = 9$, $E/K$ has potentially supersingular reduction, and $\#\bF_K \equiv 3 \bmod 4$. To prove the $p$-isogeny conjecture \eqref{main} for $E/K$ it suffices to prove it assuming that $K(E[\phi]) = K$ (without losing other assumptions).\elem

\bpf Consider the subfield $L$ of $K_\psi/K$ fixed by the $2$-Sylow subgroup of $\Gal(K_\psi/K)$. As $\Gal(K_\psi/K)$ is cyclic of order dividing $p - 1$, the degree $[K_\psi : L]$ is at most $2$ and $[L : K]$ is odd. Using \Cref{7A}, we replace $K$ by $L$ (we do not lose any assumptions by doing this; in particular,~\Cref{HH}~\ref{HHa} shows that $E/L$ still has additive reduction of Kodaira type III or III$^*$). If $K_\psi = L$ we are done, so assume that $[K_\psi : K] = 2$. The $p$-isogeny conjecture is already known for $E/K_\psi$: if $K_\psi/K$ is unramified this is \S\ref{UU}, and if it is ramified this follows from \Cref{HH}~\ref{HHa} and \Cref{six}. Using \Cref{twist} we can therefore replace $E$ by its quadratic twist $\wt{E}$ by $K_\psi/K$ (without losing any assumptions). But the Galois action on $\wt{E}[\wt{\phi}]$ is trivial by construction, so we have reduced to the case $K(E[\phi]) = K$.\epf

\bpp[Assumptions specific to the present case]\lab{S71} 
In view of \S\ref{S70} and \Cref{7B}, for the rest of the paper we will be assuming in addition to \S\ref{S61} that $v(\Delta) = 3$ or $9$, $\#\bF_K \equiv 3 \bmod 4$, the reduction is potentially good supersingular, and $K(E[\phi]) = K$. In this case the Galois closure $L$ of $F/K$ from \S\ref{S70} is of degree $8$ with $e_{L/K} = 4$, $f_{L/K} = 2$. Also, \Cref{HH}~\ref{HHa} shows that $E/L$ has good supersingular reduction. We set $G = \Gal(L/K)$ and let $I \lhd G$ be the index $2$ inertia subgroup. The subfield of $L/K$ fixed by $I$ is denoted by $M$.

The assumption $K(E[\phi]) = K$ gives in particular that $(\psi, -1)_K = 1$, so \eqref{main} in this case is
\be\lab{main2} w(E/K) \overset{?}{=} \sigma_{\phi_K}.\ee 
\epp

\bpp[A convenient Weierstrass equation] Assume the setup of \S\ref{S71} and pick a minimal Weierstrass equation for $E/K$ with associated quantities $a_1, a_2, \dotsc, c_4, c_6, \Delta, j = \f{c_4^3}{\Delta}$ (cf.~\cite{Tat75}*{\S1}). Then
\[ y^2 = x^3 - \f{c_4}{48}x - \f{c_6}{864}\]
is another minimal equation for $E/K$ since it has integral coefficients and the same valuation of the discriminant (we are assuming $p > 3$) \cite{Tat75}*{(1.7)}. If one considers it as an equation for $E/L$, it is no longer minimal but after a change of coordinates 
\[ 
\begin{split} 
x &= u^2 X, \\
y &= u^3 Y, 
\end{split}
\]
with $u = \pi_L^{v(\Delta)/3}$ one arrives at a minimal equation 
\be\lab{eqF} Y^2 = X^3 - \f{c_4}{48u^4}X - \f{c_6}{864u^6} \ee
for $E/L$. Indeed, its discriminant has valuation $0$ and it has integral coefficients because the relations $3v(c_4) \ge v(\Delta)$ (i.e.,
$v(j) \ge 0$) and $1728\Delta = c_4^3 - c_6^2$ show that $v(c_4) \ge v(\Delta)/3$ and $v(c_6) \ge v(\Delta)/2$.\epp

\bpp[The formal group of $E/L$]\lab{New} With the choice of a minimal equation \eqref{eqF} $T = -X/Y$ is a parameter for the formal group $\cF$ of $E/L$. Similarly, $t = -x/y$ is a parameter for the formal group of $E/K$, and 
\be\lab{ftr}T = ut.\ee
 Since $E/L$ has good supersingular reduction, $\cF$ is of height $2$. In other words, in $[p]_\cF(T) = pT + V_2T^2 \dotsb + V_pT^p + \dotsb + V_{p^2}T^{p^2} + \dotsb $ the first coefficient which is a unit is $V_{p^2}$. Let $\fm_0^+ \subset \ov{K}$ consist of all elements of positive valuation after uniquely extending $v$ to $\ov{K}$ (see also \S\ref{RR}). The $\Gal(\ov{K}/L)$-module $E[p]$ is isomorphic to the kernel $N_p\subset \fm_0^+$ of $[p]_\cF$ via the map $T(P) = -X(P)/Y(P)$ (one puts $T(O) = 0$). \epp

\blem\lab{odd} The $L$-valuation of a nonzero $\beta \in T(E[\phi])$ is an odd integer independent of $\beta$.\elem

\bpf  For a nonzero $P \in E[\phi]$, \eqref{ftr} shows that $v_L(T(P)) = v_L(t(P)) + v_L(u) = 4\cdot v(t(P)) + \f{v(\Delta)}{3}$ is an odd integer, because $v(t(P))$ is an integer due to our assumption that $E[\phi] \subset E(K)$. But since $E[\phi]$ and $T(E[\phi])$ are cyclic of order $p$, and the formal group law is $T_1 + T_2 + \text{higher order terms}$, all nonzero elements of $T(E[\phi])$ have the same valuation. \epf

\bpp[The actions of $G$]\lab{action} Let $\cE/\cO_L$ be the N\'{e}ron model of $E/L$, and let $\cE\pr/\cO_L$ be that of $E\pr/L$. For each $\sigma \in G = \Gal(L/K)$ we have a commutative diagram
\be\lab{neron}\begin{split} \xymatrix{
   \cE \ar[d]\ar[r]^{\sigma} &\cE \ar[d] \\
   \Spec \cO_L \ar[r]^{\sigma} &\Spec \cO_L.
}\end{split}\ee
Here $\sigma \colon \Spec \cO_L \ra \Spec \cO_L$ is a morphism corresponding to $\cO_L \xra{\sigma\i} \cO_L$ and $\sigma \colon \cE \ra \cE$ is the unique morphism making the square commute, obtained by invoking the N\'{e}ron property. Uniqueness gives us actions of $G$ on both $\Spec \cO_L$ and $\cE$ which are compatible with the morphism $\cE \ra \Spec \cO_L$. Analogous statements are true for $\cE\pr/\cO_L$.

Since $\cE/\cO_L$ is an abelian scheme, \cite{BLR90}*{\S 7.3 Proposition 6} shows that the isogeny $\phi\colon E\ra E\pr$ extends to an isogeny $\phi\colon \cE \ra \cE\pr$, whose kernel is a finite flat commutative $\cO_L$-group scheme $\cE[\phi]$ of order $p$ with generic fiber $\cE[\phi]_L = E[\phi]_L$. The diagram \eqref{neron} being functorial, we get an action of $G$ on $\cE[\phi]$ which is compatible with its action on $\Spec \cO_L$. Restricting this action to $I$ and reducing to the special fiber $\cE[\phi]_{\bF_L}$, we get an action of $I$ on $\cE[\phi]_{\bF_L}$ preserving the morphism to $\Spec \bF_L$.

Following \cite{CFKS10}*{Remark after Lemma 2.20} we define the $\cO_L/p\cO_L$-module (or $\cO_L$-module)
\be\lab{Lie} \Lie(\cE[\phi]) \ce \Ker(\cE[\phi]((\cO_L/p\cO_L)[\eps]/(\eps^2)) \ra \cE[\phi](\cO_L/p\cO_L)). \ee
(One could more accurately call this $\Lie(\cE[\phi]_{\cO_L/p \cO_L})$.) The action of $G$ on $\cE[\phi]$ gives an $\cO_L$-semilinear action of $G$ on $\Lie(\cE[\phi])$.

Since $E/L$ is a base change of $E/K$, one also has an action of $G$ on $E(L) \cong \cE(\cO_L)$ for which $\phi_L$ is $G$-equivariant, being defined over $K$. Let $\phi_{\cO_L}\colon \cE(\cO_L) \ra \cE\pr(\cO_L)$ be the map induced by $\phi\colon \cE \ra \cE\pr$ on $\cO_L$-points; it is $G$-equivariant as well. Since $E(L)^G = E(K)$, and similarly for $E\pr$, we get 
\be\lab{phi}\chi(\phi_K) = \chi(\cE(\cO_L)^G \xra{\phi_{\cO_L}} \cE\pr(\cO_L)^G).\ee
\epp

\bthm[\cite{TO70}*{pp.~14--16, Remarks 1 and 5}]\lab{OoTa} Let $A$ be $\cO_L$, $L$, $\bF_L$, or $\ov{\bF}_L$. There is a bijective correspondence between isomorphism classes of finite flat group schemes $G$ of order $p$ over $A$ and equivalence classes of factorizations $p = ac$ with $a, c \in A$, where $p = ac$ and $p = a\pr c\pr$ are said to be equivalent if there is a $u\in A^\times$ such that $a\pr = u^{p - 1}a$ and $c\pr = u^{1 - p}c $. As an $A$-scheme, the group scheme corresponding to $p = ac$ is isomorphic to $\Spec A[s]/(s^p - as)$ ($c$ appears in the description of the group law). \ethm

\brem\lab{RC} \Cref{OoTa} is part of a more general Oort--Tate classification of finite flat group schemes of order $p$, cf.~\cite{TO70}. The version stated here will be sufficient for our purposes. In \Cref{OoTa} the factorization corresponding to the constant group scheme $\underline{\bZ/p\bZ}$ is $p = 1 \cdot p$ \cite{TO70}*{pp.~8--10 and Remarks on pp.~14--15}. (With our choices for $A$, this can also be seen from \Cref{OoTa} directly, because if $\Spec A[s]/(s^p - as)$ has a nontrivial $A$-point, then $a = u^{p - 1}$ for some $u \in A$, $u \neq 0$.)\erem

\blem\lab{alphap} As finite $\bF_L$-group schemes, $\cE[\phi]_{\bF_L} \cong \gA_p$. Its corresponding factorization is $p = 0 \cdot 0$. \elem

\bpf The kernel of a $p$-isogeny between supersingular elliptic curves in characteristic $p$ is local-local. Therefore, $\cE[\phi]_{\bF_L} \cong \gA_p$, because $\gA_p$ has no twists and is the only local-local group scheme of order $p$ over $\ov{\bF}_L$. Also, $\gA_p$ is isomorphic to its own Cartier dual, so the second claim follows, because by \cite{TO70}*{p.~15, Remark~2} in characteristic $p$ Cartier duality has the effect $ac \leftrightarrow (-c)(-a)$.\epf

\bpp[The kernel of $\phi\colon \cE \ra \cE\pr$ as a scheme] Let $p = ac$ with $a, c \in \cO_L$ be a factorization corresponding to $\cE[\phi]$. Since $\cE[\phi]_L \cong E[\phi]_L$ is the constant group scheme, by \Cref{RC} its corresponding factorization is $p = 1 \cdot c_0$. \Cref{OoTa} therefore gives $a = a_0^{p - 1}$ for some $a_0 \in \cO_L$. Moreover, by \Cref{alphap} the factorization of $\cE[\phi]_{\bF_L}$ is $p = 0 \cdot 0$. We conclude that $a_0$ and $c$ are of positive valuation, and as a scheme $\cE[\phi]$ is isomorphic to $\Spec \cO_L[s]/(s^p - a_0^{p - 1}s)$ with $0 < v_L(a_0) < v_L(p)/(p - 1)$. \epp

\blem\lab{len1} We have $\length_{\cO_L} \Lie(\cE[\phi]) = (p - 1)v_L(a_0)$.\elem

\bpf Interpreting \eqref{Lie} on rings, $\Lie(\cE[\phi])$ consists of $\cO_L$-algebra homomorphisms
\[ \cO_L[s]/(s^p - a_0^{p - 1}s) \ra (\cO_L/p\cO_L)[\eps]/(\eps^2)\]
whose composite with 
\[\ba (\cO_L/p\cO_L)[\eps]/(\eps^2) &\ra \cO_L/p \cO_L, \\ \eps &\mapsto 0\ea\] 
sends $s$ to $0$. Such are given by $s \mapsto b \eps$ with $a_0^{p - 1}b = 0$, or equivalently with 
\[ b \in \pi_L^{v_L(p) - (p - 1)v_L(a_0)} \cO_L / \pi_F^{v_L(p)} \cO_L.\qedhere\] \epf

\blem\lab{len2} $v_L(a_0) = v_L(\beta)$ for any nonzero element $\beta \in T(E[\phi])$ (cf.~\Cref{odd}). \elem

\bpf One way to define the filtration $E_1(L) \supset  \dotsb \supset  E_m(L) \supset \dotsb $ discussed in \S\ref{R} is as follows (cf.~\cite{LS10}*{Lemma 5.1}): for $z \in E_1(L)$ let $S_z = \ov{\{ z\}}$ be the closure of $z$ in $\cE$; if $z\neq 0$ then $S_0 \cap S_z$ is a local Artin scheme, whose length we denote by $l(z)$; now let $E_m(L)$ consist of all $z\in E_1(L)$ with $l(z)\ge m$ (one sets $l(0) = \infty$).

One nonzero $L$-point of $\cE[\phi]$ is $s\mapsto a_0$, its closure in $\cE[\phi]$ (and hence $\cE$) is $\Spec \cO_L[s]/(s - a_0)$, the intersection with the zero section is $\Spec \cO_L[s]/(s, s - a_0)$, and the length of this local Artin scheme is $v_L(a_0)$. On the other hand, every nonzero point of $E[\phi]_L$ belongs to the filtration level $v_L(\beta)$ by definition. \epf

\blem\lab{GalD} If $N$ is a finite length $\cO_M$-module equipped with an $\cO_M$-semilinear action of $G/I$, then $\length_{\cO_K} N^{G/I} = \length_{\cO_M} N$. \elem

\bpf This is clear if $N = \bF_M \times \dotsb \times \bF_M$ by classical Galois descent for vector spaces. The general case follows by induction on the number of nonzero terms in the $(G/I)$-stable filtration $N \supset \pi_K N \supset \pi_K^2 N \supset \dotsb$ using Hilbert's theorem 90. \epf

\bcor\lab{len10} One has $\length_{\cO_K} \Lie(\cE[\phi])^G = \length_{\cO_M} \Lie(\cE[\phi])^I$. \ecor

\bpf Both lengths are finite by \Cref{len1}, and one applies \Cref{GalD} with $N = \Lie(\cE[\phi])^I$. \epf

\blem\lab{L1} $\ord_p \chi(\phi_K) \equiv \ord_p \# (\Lie(\cE[\phi])^G) \equiv \length_{\cO_M} \Lie(\cE[\phi])^I \bmod 2$. \elem

\bpf 
The second congruence holds because $f_{K/\bQ_p}$ is assumed to be odd: indeed, by \Cref{len10},
\[ 
\ord_p \# (\Lie(\cE[\phi])^G) = f_{K/\bQ_p}{\length_{\cO_K} \Lie(\cE[\phi])^G} \equiv \length_{\cO_M} \Lie(\cE[\phi])^I \bmod 2. 
\]
The first congruence is \cite{CFKS10}*{Lemma 2.20 (4) and (5)} together with \eqref{phi}. The proof given there does not use the assumption (iii) of \cite{CFKS10}*{Theorem 2.1} and therefore extends to the situation considered here. We recall the argument of op.~cit.~below.

Let $\cE_{\bF_L}$ and $\cE_{\bF_L}\pr$ denote the reductions of $E/L$ and $E\pr/L$.

\begin{claim} \lab{Claim2}
$\chi(\phi_K) = \f{\# \cE\pr(\bF_L)^G}{ \# \cE(\bF_L)^G} \cdot \f{\#(\Lie \cE_{\bF_L})^G}{\#(\Lie \cE\pr_{\bF_L})^G} \cdot \chi\p{\Lie(\cE)^G \xra{\Lie(\phi)} \Lie(\cE\pr)^G} $.
\end{claim}

\bpf
Let $\fm_L \subset \cO_L$ be the maximal ideal and choose a large $n \in \bZ_{> 0}$ such that the $G$-equivariant
\[ 
\fm_L^n\Lie \cE \ra \Ker(\cE(\cO_L) \xra{r} \cE(\cO_L/\fm_L^n))\quad  \text{and} \quad \fm_L^n\Lie \cE\pr \ra \Ker(\cE\pr(\cO_L) \xra{r\pr} \cE\pr(\cO_L/\fm_L^n)).
\]
induced by the exponential maps of $\cE$ and $\cE\pr$ are isomorphisms. By Hensel's lemma, $r$ and $r\pr$ are surjective, so \cite{Ser67}*{\S1.2, Lemma 3} and the coprimality of $\#G$ and $p$ ensure the exactness of
\[
\xymatrix{
0 \ar[r] &(\fm_L^n\Lie \cE)^G \ar[r]\ar[d] &\cE(\cO_L)^G \ar[d]^{\phi_{\cO_L}}\ar[r]^-{r} &\cE(\cO_L/\fm_L^n)^G \ar[r]\ar[d] & 0 \\
0 \ar[r] & (\fm_L^n\Lie \cE\pr)^G \ar[r] &\cE\pr(\cO_L)^G \ar[r]^-{r\pr} &\cE\pr(\cO_L/\fm_L^n)^G \ar[r] & 0 
}
\]
in spite of the presence of $G$-invariants. Therefore, \eqref{phi} together with \Cref{M,O} gives
\[
\chi(\phi_K) = \f{\# \cE\pr(\cO_L/\fm_L^n)^G}{ \# \cE(\cO_L/\fm_L^n)^G} \cdot \f{\#(\Lie \cE/\fm_L^n\Lie \cE)^G}{\#(\Lie \cE\pr/\fm_L^n\Lie \cE\pr)^G} \cdot \chi\p{\Lie(\cE)^G \xra{\Lie(\phi)} \Lie(\cE\pr)^G}.
\]
To conclude it remains to argue that one has $G$-equivariant isomorphisms 
\[
\Ker\p{\cE(\cO_L/\fm_L^{i + 1}) \ra \cE(\cO_L/\fm_L^{i})} \cong \f{\fm_L^i \Lie \cE}{\fm_L^{i + 1} \Lie \cE} \quad\quad \text{for $i \ge 1$, and similarly for $\cE\pr$.}
\]
These are supplied by deformation theory: e.g., invoke \cite{Ill05}*{Thm.~8.5.9 (a) and (the analogue of) Rem.~8.5.10 (b)} and use the zero lift to get a canonical trivialization of the appearing torsor.
\epf

\begin{claim} \lab{Claim3}
$\#(\Lie \cE_{\bF_L})^G = \#(\Lie \cE\pr_{\bF_L})^G$.
\end{claim}

\bpf
Let $D$ be the covariant Dieudonn\'{e} module of the $p$-divisible group of $\cE_{\bF_L}$, let $V$ be the Verschiebung operator of $D$, and let $D\pr$ and $V\pr$ be the corresponding objects for $\cE\pr_{\bF_L}$. The $G$-equivariant isomorphism $\Lie \cE_{\bF_L} \cong D/VD$ and \cite{Ser67}*{\S1.2,~Lemma~3} give $(\Lie \cE_{\bF_L})^G \cong D^G/VD^G$. Consequently, since $D^G$ is a free $\bZ_p$-module of finite rank and $V$ is $\bZ_p$-linear, 
\be\lab{DDVV}
\#(\Lie \cE_{\bF_L})^G = \det(V \tensor_{\bZ_p} \bQ_p\colon D^G \tensor_{\bZ_p} \bQ_p \ra D^G \tensor_{\bZ_p} \bQ_p).
\ee
It remains to note that $D^G \tensor_{\bZ_p} \bQ_p \cong (D \tensor_{\bZ_p} \bQ_p)^G$, so the right hand side of \eqref{DDVV} is the same for $\cE\pr_{\bF_L}$ because $\phi$ induces a $G$-isomorphism $(D \tensor_{\bZ_p} \bQ_p, V \tensor_{\bZ_p} \bQ_p) \cong (D\pr \tensor_{\bZ_p} \bQ_p, V\pr \tensor_{\bZ_p} \bQ_p)$.
\epf

\begin{claim} \lab{Claim4}
$\ord_p \# \cE(\bF_L)^G = \ord_p \# \cE\pr(\bF_L)^G$.
\end{claim}

\bpf
Since the action of the inertia $I \lhd G$ preserves the morphism $\cE_{\bF_L} \ra \Spec \bF_L$, the subgroup $(\cE(\ov{\bF}_L)[p\I])^I \subset \cE(\ov{\bF}_L)[p\I]$ is $\Gal(\ov{\bF}_L/\bF_L)$-stable. Let $\Frob_K \in \Gal(\ov{K}/K)$ be a geometric Frobenius. On the one hand, the action of $\Frob_K$ on $(\cE(\ov{\bF}_L)[p\I])^I$ lifts the action of the generator of $G/I$. On the other, the action of $\Frob_K^2$ is that of the geometric Frobenius in $\Gal(\ov{\bF}_L/\bF_L)$. In conclusion, 
\[
\cE(\bF_L)^G[p\I] = \Ker\p{1 - \Frob_K\colon (\cE(\ov{\bF}_L)[p\I])^I \ra (\cE(\ov{\bF}_L)[p\I])^I}.
\]
Since $(\cE(\ov{\bF}_L)[p\I])^I \subset \cE(\ov{\bF}_L)[p\I]$ is cut out by an idempotent of $\bZ_p[I]$, it inherits $p$-divisibility. Set $T_p \ce \varprojlim\, (\cE(\ov{\bF}_L)[p^n])^I$ and $V_p \ce T_p \tensor_{\bZ_p} \bQ_p$. Since $\cE(\bF_L)^G$ is finite, $1 - \Frob_K$ is injective on $T_p$, and hence also on $V_p$. Consequently, snake lemma applied to 
\[
\xymatrix{
0 \ar[r]  & T_p \ar[r]\ar[d]_{1 - \Frob_K} & V_p \ar[r]\ar[d]^{1 - \Frob_K} & (\cE(\ov{\bF}_L)[p\I])^I \ar[r]\ar[d]^{1 - \Frob_K} & 0 \\
0 \ar[r]  & T_p \ar[r] & V_p \ar[r] & (\cE(\ov{\bF}_L)[p\I])^I \ar[r] & 0
}
\]
gives the first equality in 
\be\lab{TpVp}
\ord_p \# \cE(\bF_L)^G = \ord_p \#\p{\f{T_p}{(1 - \Frob_K)T_p}} = \ord_p \det(1 - \Frob_K\colon V_p \ra V_p).
\ee
The claim follows from \eqref{TpVp}: indeed, similar reasoning applies to $\cE\pr$, so, denoting by $V_p\pr$ the analogue of $V_p$, one notes that $\phi$ induces a $\Frob_K$-equivariant isomorphism $V_p \cong V_p\pr$.
\epf

\begin{claim} \lab{Claim1}
$\chi\p{\Lie(\cE)^G \xra{\Lie(\phi)} \Lie(\cE\pr)^G} = \#\p{\Lie(\cE[\phi]_{\cO_L/p\cO_L})^G} \overset{\S\ref{action}}{=} \#(\Lie(\cE[\phi])^G)$.
\end{claim}

\bpf
The Lie algebras $\Lie(\cE)$ and $\Lie(\cE\pr)$ are free $\cO_L$-modules of rank $1$. Consideration of the isogeny dual to $\phi$ shows that $\Lie(\cE) \xra{\Lie(\phi)} \Lie(\cE\pr)$ is injective and its cokernel $Q$ is killed by $p$. 

Consider the short exact sequence
\be\lab{mod-p}
0 \ra \cE[\phi]_{\cO_L/p\cO_L} \ra \cE_{\cO_L/p\cO_L} \ra \cE\pr_{\cO_L/p\cO_L} \ra 0
\ee
of $\cO_L/p\cO_L$-group schemes. Forming Lie algebras commutes with base change, so \eqref{mod-p} gives the~exact
\[
0 \ra \Lie(\cE[\phi]_{\cO_L/p\cO_L}) \ra \Lie(\cE)\tensor_{\cO_L} \cO_L/p\cO_L \xra{\Lie(\phi) \tensor_{\cO_L} \cO_L/p\cO_L} \Lie(\cE\pr)\tensor_{\cO_L} \cO_L/p\cO_L.
\]
Consequently, snake lemma applied to the commutative diagram
\[
\xymatrix{
0 \ar[r] & \Lie(\cE) \ar[r]^-{\Lie(\phi)}\ar[d]^{p} & \Lie(\cE\pr) \ar[r]\ar[d]^{p} & Q \ar[r]\ar[d]^{0} & 0 \\
0 \ar[r] & \Lie(\cE) \ar[r]^-{\Lie(\phi)} & \Lie(\cE\pr) \ar[r] & Q \ar[r] & 0
}
\]
of $G$-modules gives $Q \cong \Lie(\cE[\phi]_{\cO_L/p\cO_L})$. Since $(\#G, p) = 1$, the resulting
\[
0 \ra \Lie(\cE) \xra{\Lie(\phi)} \Lie(\cE\pr) \ra \Lie(\cE[\phi]_{\cO_L/p\cO_L}) \ra 0
\]
remains short exact after taking $G$-invariants, and the desired conclusion follows.
\epf

\Cref{Claim1,Claim2,Claim3,Claim4} provide an equality underlying the first congruence of \Cref{L1}.
\epf

\bpp[Vector space structures]\lab{maxnr} Let $M_0 \ce \wh{M^{\text{ur}}}$ be the completion of the maximal unramified extension of $M$, and let and $L_0$ be a compositum of $M_0$ and $L$. The field $L_0$ is complete, and $L_0/M_0$ is Galois with Galois group identified with $I$. In particular, $[L_0 : M_0] = 4$. The rings of integers of $L_0$ and $M_0$ will be denoted by $\cO_{L_0}$ and $\cO_{M_0}$. 

As observed in \S\ref{action}, $I$ acts on $\cE[\phi]$ compatibly with its action on $\Spec \cO_L$. With the identification above, $I$ therefore acts on $\cE[\phi]_{\cO_{L_0}}$ compatibly with its action on $\Spec \cO_{L_0}$, and we get an $\cO_{L_0}$-semilinear action of $I$ on $\Lie(\cE[\phi]_{\cO_{L_0}}) \cong  \Lie(\cE[\phi]) \tensor_{\cO_L} \cO_{L_0}$. 

By definition $\Lie(\cE[\phi])$ is an $\cO_L/p\cO_L$-module. In particular, denoting by $W(\bF_L)$ the ring of Witt vectors, we can regard $\Lie(\cE[\phi])$ as a vector space over $\bF_L \cong W(\bF_L)/pW(\bF_L) \subset \cO_L/p\cO_L$ equipped with an $\bF_L$-linear action of $I$. In other words, $\Lie(\cE[\phi])$ is a finite dimensional $\bF_L$-representation of $I$ with 
\be\lab{Aaa} \length_{\cO_L} \Lie(\cE[\phi]) = \dim_{\bF_L} \Lie(\cE[\phi]),\ee
and also
\be\lab{Bbb} \length_{\cO_M} \Lie(\cE[\phi])^I = \dim_{\bF_L} \Lie(\cE[\phi])^I.\ee
On the other hand, $\Lie(\cE[\phi]_{\cO_{L_0}}) \cong \Lie(\cE[\phi]) \tensor_{\bF_L} \ov{\bF}_L$, and also $\Lie(\cE[\phi]_{\cO_{L_0}})^I \cong \Lie(\cE[\phi])^I \tensor_{\bF_L} \ov{\bF}_L$; therefore
 \be\lab{Ccc} \length_{\cO_L} \Lie(\cE[\phi]) \overset{\eqref{Aaa}}{=} \dim_{\bF_L} \Lie(\cE[\phi]) = \dim_{\ov{\bF}_L} \Lie(\cE[\phi]_{\cO_{L_0}}),\ee
and also
\be\lab{len3} \length_{\cO_M} \Lie(\cE[\phi])^I \overset{\eqref{Bbb}}{=} \dim_{\bF_L} \Lie(\cE[\phi])^I = \dim_{\ov{\bF}_L} \Lie(\cE[\phi]_{\cO_{L_0}})^I.\ee \epp

\bcor\lab{I1} $\ord_p \chi(\phi_K) \equiv \dim_{\ov{\bF}_L} \Lie(\cE[\phi]_{\cO_{L_0}})^I \bmod 2$.\ecor

\bpf Combine \Cref{L1} and \eqref{len3}.\epf

\bpp[The Dieudonn\'{e} module of the special fiber]\lab{DDN}The special fiber of $\cE[\phi]_{\cO_{L_0}}$ is $\cE[\phi]_{\ov{\bF}_L}$, which by \S\ref{action} carries the action of $I$ preserving the morphism to $\Spec \ov{\bF}_L$. By \Cref{alphap}, $\cE[\phi]_{\ov{\bF}_L} \cong \gA_p$, so the (covariant) Dieudonn\'{e} module $D(\cE[\phi]_{\ov{\bF}_L}) \cong D(\gA_p)$ is especially easy to describe: $D(\cE[\phi]_{\ov{\bF}_L}) \cong~\ov{\bF}_L$ with vanishing Frobenius and Verschiebung.  By functoriality, $D(\cE[\phi]_{\ov{\bF}_L})$ is an $\ov{\bF}_L$-representation of $I$. The latter is cyclic of order $4$, so it acts on $D(\cE[\phi]_{\ov{\bF}_L})$ via scaling by some $4^{\mathrm{th}}$ roots of unity.\epp

\bpp[$\ov{\bF}_M$-representations of inertia] \lab{RR} Let $I_M \lhd \Gal(\ov{K}/M)$ be the inertia subgroup, and let $P_M \lhd I_M$ be the wild inertia. We are interested in continuous irreducible $\ov{\bF}_M$-representations $V$ of $I_M$. Since $\Char \ov{\bF}_M = p$, and $P_M$ is pro-$p$, one has $V^{P_M} \neq 0$ \cite{Ser77}*{Proposition 26}. Moreover, $P_M$ is normal in $I_M$, so $V^{P_M}$ is $I_M$-stable, hence $V^{P_M} = V$. In other words, $V$ is the inflation of a continuous irreducible representation of the tame inertia $I_M/P_M$. Tame inertia is abelian, so $V$ is $1$-dimensional; it must, therefore, be isomorphic to some $V_a$, $a \in \bQ$ constructed as follows (cf.~\cite{Ser72}*{\S\S1.7--1.8}). The valuation $v$ on $K$ extends uniquely to a (no longer discrete) valuation on $\ov{K}$, which we continue to denote $v$. Let $\fm_a$ be the set of $x\in \ov{K}$ with valuation $v(x)\ge a$. Let $\fm_a^+ \subset \fm_a$ be the set of $x$ with $v(x) > a$. The quotient $\fm_a/\fm_a^+$ is a $1$-dimensional $\fm_0/\fm_0^+ \cong \ov{\bF}_M$-linear representation of $I_M$, which we call $V_a$. In addition, $V_a \cong V_b$ if and only if $a - b \in \bZ[1/p]$ (loc.~cit.), and we conclude that the Grothendieck group of continuous $\ov{\bF}_M$-representations of $I_M$ is isomorphic to the group ring $R \ce \bZ[\bQ/\bZ[1/p]]$. Multiplication in $R$ corresponds to tensor product of representations (loc.~cit.).

In fact, since $I_M$ identifies with $\Gal(\ov{M_0}/M_0)$, we can think of $R$ as the Grothendieck group of continuous $\ov{\bF}_M$-linear representations of $\Gal(\ov{M_0}/M_0)$. The representations that will interest us most are $\Lie(\cE[\phi]_{\cO_{L_0}})$ and $D(\cE[\phi]_{\ov{\bF}_L})$; they factor through the finite quotient $\Gal(\ov{M_0}/M_0)/\Gal(\ov{M_0}/L_0) = I$ of order $4$.\epp

\bpp[Maps related to $R$]\lab{maps} Following \cite{CFKS10}*{\S 7.2} let us supplement the ring $R = \bZ[\bQ/\bZ[1/p]]$ of \S\ref{RR} with

\begin{enumerate}[label={(\alph*)}]
\item The notation $\gamma(a)$ for the standard $\bZ$-basis element of $R$ corresponding to $a\in \bQ/\bZ[1/p]$;
\item The automorphism $\varphi\colon R \ra R$ induced by sending $\gamma(a)$ to $\gamma(pa)$;
\item The $\bZ$-linear map $\gA\colon R \ra \bQ/\bZ[1/p]$ which sends $\gamma(a)$ to $a$;
\item The $\bZ$-linear map $\wt{\gA}\colon R \ra \bQ$ defined by sending $\gamma(a)$ to the unique element of $\bZ_{(p)}\cap (0, 1]$ whose class $\bmod\ \bZ[1/p]$ is $a$;
\item The $\bZ$-linear map $\delta_0\colon R \ra \bZ$ such that $\delta_0(\gamma(0)) = 1$ and $\delta_0(\gamma(a)) = 0$ for $a \neq 0$;
\item The $\bZ$-linear map $\deg\colon R \ra \bZ$ which sends each $\gamma(a)$ to $1$.
\end{enumerate}
Thinking in terms of representations of $\Gal(\ov{M_0}/M_0)$, one observes that $\deg$ (resp.,~$\delta_0$) is nothing else than the $\ov{\bF}_M$-dimension of the representation space (resp.,~the fixed subspace).\epp

\bprop\lab{App} Denoting by $[V]$ the class in $R$ of an $\ov{\bF}_M$-representation $V$ of $\Gal(\ov{M_0}/M_0)$, we have the following relations 
\begin{enumerate}[label={(\alph*)}]
\item\lab{Appa} $\displaystyle{ \wt{\gA}(\varphi\i([D(\cE[\phi]_{\ov{\bF}_L})])) - \wt{\gA}([D(\cE[\phi]_{\ov{\bF}_L})]) = \f{\deg([\Lie(\cE[\phi]_{\cO_{L_0}})])}{4} - \delta_0([\Lie(\cE[\phi]_{\cO_{L_0}})]) }$ in $\bQ$;
\item\lab{Appb} $\displaystyle{-\gA([D(\cE[\phi]_{\ov{\bF}_{L}})]) = \f{p\cdot  \deg([\Lie(\cE[\phi]_{\cO_{L_0}})])}{4(p - 1)}}$ in $\bQ/\bZ[1/p]$.
\end{enumerate} \eprop

\bpf This is \cite{CFKS10}*{Proposition 7.3} applied to $P = \cE[\phi]_{\cO_{L_0}}$; their $K$ is our $M_0$, their $L$ is our $L_0$, their $k$ is our $\ov{\bF}_M$ ($= \ov{\bF}_L$), and their $\Delta$ is our $I$. Since $E[\phi] \subset E(K)$, the Galois representation on geometric points of the generic fiber of $\cE[\phi]_{\cO_{L_0}}$ is trivial, which allows us to discard the first summand in \cite{CFKS10}*{Proposition 7.3 (2)}.  \epf

\bprop\lab{three} Under the assumptions of \S\ref{S71}, the $p$-isogeny conjecture \eqref{main2} is true. \eprop

\bpf From \S\ref{DDN} we get that $[D(\cE[\phi]_{\ov{\bF}_{L}})]$ is $\gamma(i)$, where $i = 0, \f{1}{4}, \f{1}{2}$, or $\f{3}{4}$. Moreover, 
\[ \deg([\Lie(\cE[\phi]_{\cO_{L_0}})]) \overset{\S\ref{maps}}{=} \dim_{\ov{\bF}_L} \Lie(\cE[\phi]_{\cO_{L_0}}) \overset{\eqref{Ccc}}{=} \length_{\cO_L} \Lie(\cE[\phi]) \overset{\ref{len1}}{=} (p - 1)v_L(a_0) =  (p - 1)(2m + 1),\] 
for some $m \ge 0$, where the last equality follows from \Cref{len2} and \Cref{odd}. \Cref{App}~\ref{Appb} gives 
\[ -i = \f{p(2m + 1)}{4} \quad \text{in } \bQ/\bZ[1/p]. \] 
Since $p \equiv 3 \bmod 4$ this means that $i = \f{1}{4}$ if $m$ even, and $i = \f{3}{4}$ if $m$ is odd. Therefore, $\varphi\i(\gamma(i)) = \gamma(\f{3}{4})$ if $m$ is even, and $\varphi\i(\gamma(i))= \gamma(\f{1}{4})$ if $m$ is odd. The left hand side of \Cref{App}~\ref{Appa} is therefore $\f{1}{2}$ if $m$ is even, and $-\f{1}{2}$ if $m$ is odd.

Write $p = 4k + 3$. If $m$ is even, \Cref{App}~\ref{Appa} gives 
\[  \delta_0([\Lie(\cE[\phi]_{\cO_{L_0}})]) = \f{(p - 1)(2m + 1)}{4} - \f{1}{2} = k(2m + 1) + m \equiv k \bmod 2.\]
If $m$ is odd, it gives
\[  \delta_0([\Lie(\cE[\phi]_{\cO_{L_0}})]) = \f{(p - 1)(2m + 1)}{4} + \f{1}{2} = k(2m + 1) + (m + 1) \equiv k \bmod 2.\]
We conclude that in all cases $\delta_0([\Lie(\cE[\phi]_{\cO_{L_0}})]) \equiv k \bmod 2$. On the other hand,
\[ \delta_0([\Lie(\cE[\phi]_{\cO_{L_0}})]) \overset{\S\ref{maps}}{=} \dim_{\ov{\bF}_L} \Lie(\cE[\phi]_{\cO_{L_0}})^I \overset{\ref{I1}}{\equiv} \ord_p \chi(\phi_K) \bmod 2,  \] 
so $\sigma_{\phi_K} = (-1)^k$. 

To compute the root number, note that by \Cref{RN}~\ref{RNe}, $w(E/K) {=} (-1)^{\lfloor v(\Delta) \cdot \# \bF_K/12 \rfloor}$. The latter is $(-1)^{\lfloor v(\Delta) \cdot p/12 \rfloor}$, because $\# \bF_K = p^{f_{K/\bQ_p}} \equiv p \bmod 24$, since $f_{K/\bQ_p}$ is odd. Because $v(\Delta) \in \{3, 9\}$, one checks that $w(E/K) = 1$ if $p \equiv 3 \bmod 8$, and $w(E/K) = -1$ if $p \equiv 7 \bmod 8$. The former occurs if $k$ is even, the latter if $k$ is odd. \epf

\bthm\lab{DONE} The $p$-isogeny conjecture is true for $p > 3$.\ethm

\bpf Indeed, we have settled all the remaining cases; see \S\ref{sum}, \Cref{six,two,four}, \S\ref{S70}, \Cref{7B}, \S\ref{S71}, and \Cref{three}. \epf

\section{The $p$-parity conjecture for elliptic curves with complex multiplication}\lab{S7}

In this section $E$ denotes an elliptic curve over a number field $K$ such that $E/K$ has complex multiplication by an order of the imaginary quadratic field $F \ce (\End_K E) \tensor \bQ$. The ring of integers of $F$ will be denoted by $\cO_F$. We set $\cX_p(E/K) \ce \Hom_{\bZ_p} (\varinjlim \Sel_{p^n}(E/K), \bQ_p/\bZ_p) \tensor_{\bZ_p} \bQ_p$ and note that by definition $\rk_p(E/K) = \dim_{\bQ_p} \cX_p(E/K)$.

In \Cref{CM-A} we prove the $p$-parity conjecture for such elliptic curves; I thank Karl Rubin for pointing out to me that this result follows from \Cref{G}. In \Cref{CM-B} we prove that the global root number of $E/K$ is $1$, which gives \Cref{CM-Ba}. We begin by recalling two well-known results that will be used in the proofs.

\bprop\lab{CM-iso} There is an elliptic curve $E\pr/K$ with $\End_K E\pr \cong \cO_F$ and an isogeny $\gL \colon E \ra E\pr$ defined over $K$. \eprop

\bpf See, for instance, \cite{Rub99}*{Proposition 5.3}. \epf

\bprop\lab{CM-iso2} If $E$ and $E\pr$ are two elliptic curves defined over a number field $K$ and $\gL \colon E \ra E\pr$ is an isogeny defined over $K$, then $\rk_p(E/K) = \rk_p(E\pr/K)$ and $w(E/K) = w(E\pr/K)$. \eprop

\bpf If $\gL\pr$ is the dual isogeny, one notes that $\gL\pr \circ \gL$ induces automorphisms of $\cX_p(E/K)$ and $V_l(E/K)$, and similarly for $\gL \circ \gL\pr$. Hence the maps induced by $\gL$ are isomorphisms. This gives the claim about $\rk_p$.  The local root numbers at finite places are defined in terms of $V_l(E/K)$ and at infinite places are $-1$ by \Cref{RN}~\ref{RNa}, so the conclusion follows. \epf

\bprop\lab{CM-B} If $E$ has complex multiplication defined over $K$, then $w(E/K) = 1$. More precisely, if $\psi_{E/K} = \prod\pr_v \psi_v \colon \bA_K^\times/K^\times \ra \bC^\times$ is the Hecke character associated to $E/K$ (cf., for instance, \cite{Rub99}*{Theorem 5.15}), then $w(E/K_v) = \psi_v(-1)$ for every place $v$. \eprop

\bpf It is clear that the first claim follows from the second by taking the product over all places: 
\[ \prod_v \psi_v(-1) = \psi_{E/K}(-1) = 1.\]
Let $v$ be a finite place of $K$ and choose a rational prime $p$ such that $v \nmid p$ and $p$ splits in $F$. The Galois representation $V_p(E/K_{v})$ is a direct sum $\psi_{v} \oplus \psi_{v}$ \cite{Rub99}*{Corollary 5.6 and Theorem 5.15~(ii)} (here we engage in the usual abuse of local class field theory by identifying characters of $K_{v}^\times$ and of $\cW(\ov{K}_{v}/K_{v})$). If $\omega_{v}\colon \cW(\ov{K}_{v}/K_{v}) \ra \bC^\times$ is the cyclotomic character, then, because of the Weil pairing, $\psi_{v}\omega_{v}^{-1/2}$ squares to the trivial character, and hence is self-contragredient. The determinant formula \ref{EEc} from \S\ref{EE} applies, giving $w(\psi_{v}\omega_{v}^{-1/2}, \eta)^2 = \psi_{v}(-1)$. Since a twist by $\omega_{v}^{-1/2}$ does not affect the root number \cite{Roh94}*{\S11 Proposition (iii)}, we conclude that $w(E/K_{v}) = w(\psi_{v}, \eta)^2 = \psi_{v}(-1)$.

The formula $w(E/K_v) = \psi_v(-1)$ holds at an archimedean place $v$ as well. Indeed, $w(E/K_v) = -1$ by \Cref{RN}~\ref{RNa}, while $\psi_v(-1) = -1$ by construction of $\psi_{E/K}$, see the proof of \cite{Rub99}*{Theorem~5.15} (the $\psi_{E/K}$ constructed there is unique because (ii) there determines its finite component uniquely, and then the infinite component is uniquely determined because $\psi_{E/K}$ is a Hecke character).  \epf

\bthm\lab{CM-A} If $E$ has complex multiplication defined over $K$, then the $p$-parity conjecture holds for $E/K$. \ethm

\bpf Due to \Cref{CM-iso,CM-iso2}, we assume that $E$ has complex multiplication by the maximal order $\cO_F$. If $p$ is inert or ramifies in $F$, then $F_{\fp} \ce \cO_F \tensor \bQ_p$ is a quadratic extension of $\bQ_p$. Since $\cX_p(E/K)$ is an $F_{\fp}$-vector space, $\rk_p(E/K) = \dim_{\bQ_p} \cX_p(E/K)$ is even and the conclusion follows from \Cref{CM-B}. On the other hand, if $p$ ramifies or splits in $F$ and $\fp$ is a prime of $F$ above $p$, then $E[\fp] \ce \bigcap_{\gA \in \fp} \{ P \in E(\ov{K})\colon \gA P = 0\}$ is a subgroup of $E[p]$ of order $p$ defined over $K$ (see \cite{Rub99}*{Proposition 5.4} for instance). In other words, $E[\fp]$ is the kernel of a $p$-isogeny defined over $K$, and the conclusion follows from \Cref{G}. \epf

\end{document}